\documentclass[reqno,10pt]{amsart}
\usepackage{amsmath,amsthm,amssymb}
\usepackage[dvipdfmx]{graphicx}
\usepackage{enumerate}
\usepackage{color}
\usepackage{mathrsfs}
\newtheorem{Th}{Theorem}[section] 
\newtheorem{Lem}[Th]{Lemma} 
\newtheorem{Prop}[Th]{Proposition}
\newtheorem{Cor}[Th]{Corollary} 
 
\newtheorem{Prob}[Th]{Problem} 
\newtheorem{Rem}[Th]{Remark}  
\newtheorem{Hyp}[Th]{Hypothesis} 

\renewcommand{\theequation}{\arabic{section}.\arabic{equation}}

\newcommand{\R}{{\mathbb R}}
\newcommand{\N}{{\mathbb N}}
\newcommand{\Hb}{{H^1(\Omega)}}

\newcommand{\Hsecond}{{H^1_{\Gamma_2}(\Omega)}}

\newcommand{\Lo}{{L^2(\Omega)}}
\newcommand{\Cinfb}{{C^\infty(\overline{\Omega})}}
\newcommand{\Cinfo}{{C^\infty_0(\Omega)}}
\newcommand{\dvg}{\mbox{\rm div}\,}
\newcommand{\grad}{\nabla}
\newcommand{\rot}{\nabla\times}
\newcommand{\into}{\int_\Omega}

\newcommand{\eps}{\varepsilon}
\newcommand{\DS}{\displaystyle}

\newcommand{\upre}{u^*}
\newcommand{\vpre}{v^*}

\makeatletter 
\long\def\@makefntext#1{\parindent 1em\noindent 
\@hangfrom{\hbox to 1.8em{\hss$^{\@thefnmark}$}}#1}
\makeatother

\begin{document}

\title[Projection method with the total pressure boundary condition]{
    A projection method for Navier--Stokes equations 
    with\\ a boundary condition including the total pressure}

\author{Kazunori Matsui}
\address[Kazunori Matsui]{Division of Mathematical and Physical Sciences, 
Graduate School of Natural Science and Technology, Kanazawa University, 
Kanazawa 920-1192, Japan}
\email{first-lucky@stu.kanazawa-u.ac.jp}

\date{\today}

\keywords{Projection method, The Navier--Stokes equations,
Total pressure, Finite element method}

\subjclass[2020]{\emph{Primary}: 65M12; \emph{Secondary}: 35Q30, 76D03, 76D05, 76M10} 

\begin{abstract} 
    We consider a projection method for time-dependent 
    incompressible Navier--Stokes equations with 
    a total pressure boundary condition. 
    The projection method is one of the numerical calculation methods 
    for incompressible viscous fluids often used in engineering. 
    In general, the projection method needs 
    additional boundary conditions to solve a pressure-Poisson equation, 
    which does not appear in the original Navier--Stokes problem. 
    On the other hand, many mechanisms generate flow 
    by creating a pressure difference, such as 
    water distribution systems and blood circulation. 
    We propose a new additional boundary condition 
    for the projection method with a Dirichlet-type pressure 
    boundary condition and no tangent flow. 
    We demonstrate stability for the scheme and 
    establish error estimates for the velocity and pressure 
    under suitable norms. A numerical experiment verifies 
    the theoretical convergence results. 
    Furthermore, the existence of a weak solution 
    to the original Navier--Stokes problem is proved 
    by using the stability.
\end{abstract} 

\maketitle


\section{Introduction}
\label{intro}
Let $T>0$ and let $\Omega$ be a bounded Lipschitz domain in $\R^d~(d=2,3)$
with the boundary $\Gamma=\Gamma_1\cup\Gamma_2$
(see also Section \ref{sec_notation} for the precise assumption).
We consider the following Navier--Stokes problem: 
Find two functions 
$u:\Omega\times[0,T]\rightarrow\R^d$ and $p:\Omega\times[0,T]\rightarrow\R$ such that
\begin{align}\label{original}\left\{\begin{array}{ll}
    {\DS \frac{\partial u}{\partial t}+(u\cdot\nabla)u
    -\nu\Delta u+\frac{1}{\rho}\grad p=f}&\mbox{in }\Omega\times(0,T),\\[8pt]
    \dvg u=0&\mbox{in }\Omega\times(0,T),\\
    u=0&\mbox{on }\Gamma_1\times(0,T),\\
    u\times n=0&\mbox{on }\Gamma_2\times(0,T),\\
    {\DS p+\frac{\rho}{2}|u|^2=p^b}&\mbox{on }\Gamma_2\times(0,T),\\
    u(0)=u_0&\mbox{in }\Omega,
\end{array}\right.\end{align}
where $\nu,\rho>0$, $f:\Omega\times(0,T)\rightarrow\R^d$, 
$p^b:\Gamma_2\times(0,T)\rightarrow\R$,
$u_0:\Omega\rightarrow\R^d$,
$n$ is the unit outward normal vector for $\Gamma$
and ``$\times$'' is the cross product in $\R^d$
\footnote{If $d=2$, then we define 
\[
    v\times w:=v_x w_y-v_y w_x\in\R\qquad\mbox{for all }v=(v_x,v_y),w=(w_x,w_y)\in\R^2.
\]}.
The functions $u$ and $p$ are the velocity and the pressure of the flow
governed by (\ref{original}), respectively.
For $\Gamma_2$, we assume a boundary condition including 
a pressure value $p+\frac{\rho}{2}|u|^2$,
which is called the total pressure or stagnation pressure.
Usual pressure is often called static pressure to distinguish it 
from the total pressure.
In an experimental measurement of the total and static pressure
using a Pitot tube, the boss measurement
is dependent on the yaw angle of the Pitot tube.
Then, the effect on the total pressure $p+\frac{\rho}{2}|u|^2$ is smaller than
the effect on the usual pressure $p$ \cite[Section 7.15]{Holman01}.
The boundary condition on $\Gamma_2$ in (\ref{original}) 
is introduced in \cite{BCMP87}, and the existence of 
a weak velocity solution is proved in \cite{Bernard03,KC15}.
We will show the existence in a different way (Corollary \ref{cor_existsol}).
The stationary case has been studied in 
\cite{BCMP88,Bernard02,BRY15,CMP94,KC15}. 
In \cite{BRY15,BCPS17}, the finite element discretization problems 
with this type of boundary condition are proposed.

Next, we introduce a projection method for (\ref{original}).
The projection method is one of the numerical schemes
for Navier--Stokes equations \cite{Chorin68,Temam69}.
Error analysis in the case of the full Dirichlet boundary condition 
for the velocity is carried out in \cite{BC07,Prohl97,Rannacher92,Shen92,Shen94}. 
In the case of a boundary condition for the static pressure, 
the finite element analysis of a projection method is proposed
in \cite{GQ97,GQ98}. 
For the nonlinear term in the first equation of (\ref{original}),
it holds that
\[
    (u\cdot\nabla)u=(\nabla\times u)\times u+\frac{1}{2}\nabla|u|^2
\]
(cf. \cite{Gresho91}).
Hence, if we put $D(v,w):=(\nabla\times v)\times w$ and
$P=p+\frac{\rho}{2}|u|^2$,
then (\ref{original}) is equivalent to the following\footnote{
  If $d=2$, then $\nabla\times v$ and $(\nabla\times v)\times w$ 
  denote the scalar and vector functions, respectively, 
  defined as follows:
  for all $v=(v_x,v_y),w=(w_x,w_y)\in\R^2$,
  \[
      \nabla\times v
      :=\partial_x v_y-\partial_y v_x,\quad
      (\nabla\times v)\times w
      :=(w_y(\partial_y v_x-\partial_x v_y),
      w_x(\partial_x v_y-\partial_y v_x)).
  \]
}:
\begin{align}\label{original2}\left\{\begin{aligned}
    &\frac{\partial u}{\partial t}+D(u,u)
    -\nu\Delta u+\frac{1}{\rho}\grad P=f
    &&\mbox{in }\Omega\times(0,T),\\
    &\dvg u=0&&\mbox{in }\Omega\times(0,T),\\
    &u=0&&\mbox{on }\Gamma_1\times(0,T),\\
    &u\times n=0&&\mbox{on }\Gamma_2\times(0,T),\\
    &P=p^b&&\mbox{on }\Gamma_2\times(0,T),\\
    &u(0)=u_0&&\mbox{in }\Omega
\end{aligned}\right.\end{align}
In \cite{Cho16}, a projection method for the rotation form 
(the first equation of (\ref{original2}))
using the total pressure is introduced 
to avoid checkerboard oscillation of pressure 
in the finite difference method.

Let $\tau(:=T/N<1,N\in\N)$ be a time increment and 
let $t_k:=k\tau~(k=0,1,\ldots,N)$.
We set $\upre_0:=u_0$ and
calculate $\upre_k,u_k,p_k~(k=1,2,\ldots,N)$
by repeatedly solving the following problems (Step 1) and (Step 2).

(Step 1)
Find $\upre_k:\Omega\rightarrow\R^d$ such that
\begin{align}\label{step1}\left\{\begin{aligned}
    &\frac{\upre_k-u_{k-1}}{\tau}+D(\upre_{k-1},\upre_k)
    -\nu\Delta\upre_k=f(t_k)&&\mbox{in }\Omega,\\
    &\upre_k=0&&\mbox{on }\Gamma_1,\\
    &\upre_k\times n=0&&\mbox{on }\Gamma_2,\\
    &\dvg \upre_k=0&&\mbox{on }\Gamma_2.
\end{aligned}\right.\end{align}

(Step 2)
Find $P_k:\rightarrow\R$ and $u_k:\rightarrow\R^d$ such that
\begin{align}\label{step2_1}\left\{\begin{aligned}
    &-\frac{\tau}{\rho}\Delta P_k=-\dvg\upre_k&&\mbox{in }\Omega,\\
    &\frac{\partial P_k}{\partial n}=0&&\mbox{on }\Gamma_1,\\
    &P_k=p^b(t_k)&&\mbox{on }\Gamma_2,
\end{aligned}\right.\\
\begin{aligned}\label{step2_2}
    u_k=\upre_k-\frac{\tau}{\rho}\nabla P_k\quad\mbox{in }\Omega.
\end{aligned}\end{align}

For the velocity boundary condition on $\Gamma_2$,
we can rewrite the third and fourth equations of (\ref{step1}) 
by using $\kappa:=\dvg n=(d-1)\times $(mean curvature) as follows:

\begin{Rem}\label{rem_KN}
    If $v\in C^1(\overline{\Omega})$ satisfies that $v\times n=0$
    on $\Gamma_2$, then we have
    \[
        \frac{\partial v}{\partial n}\cdot n+\kappa v\cdot n
        =\dvg v\qquad\mbox{on }\Gamma_2.
    \]
    For the proof, see \cite[Lemma 7]{KN99}.
    Hence, the third and fourth equations of (\ref{step1}) 
    are equivalent to the following equations:
    \[
        \upre_k\times n=0,\quad
        \frac{\partial\upre_k}{\partial n}\cdot n+\kappa \upre_k\cdot n=0
        \qquad\mbox{on }\Gamma_2.
    \]
    In particular, if $\Gamma_2$ is flat, then it holds that
    \[
        \upre_k\times n=0,\quad
        \frac{\partial\upre_k}{\partial n}\cdot n=0
        \qquad\mbox{on }\Gamma_2.
    \]
\end{Rem}

\begin{Rem}
    By replacing $u_{k-1}$ in the first equation of (\ref{step1})
    with (\ref{step2_2}) at the previous step,
    it holds that for all $k=1,2,\ldots,N$,
    \[
        \frac{\upre_k-\upre_{k-1}}{\tau} + D(\upre_{k-1},\upre_k)
        - \nu\Delta\upre_k + \frac{1}{\rho}\nabla P_{k-1}
        = f(t_k)\qquad \mbox{in }\Omega.
    \]
    It follows from (\ref{step2_1}) and (\ref{step2_2}) that 
    $\dvg u_k=0$ in $\Omega$, $u_k\cdot n=0$ on $\Gamma_1$.
    Hence, by (\ref{step1}), (\ref{step2_1}), and (\ref{step2_2}),
    it holds that for all $k=1,2,\ldots,N$,
    \begin{align*}\left\{\begin{aligned}
        &\frac{\upre_k-\upre_{k-1}}{\tau} + D(\upre_{k-1},\upre_k)
        - \nu\Delta\upre_k + \frac{1}{\rho}\nabla P_{k-1}
        = f(t_k)&&\mbox{in }\Omega,\\
        &\dvg u_k=0 &&\mbox{in }\Omega,\\
        &\upre_k=0 &&\mbox{on }\Gamma_1,\\
        &\upre_k\times n=0&&\mbox{on }\Gamma_2,\\
        &P_k=p^b(t_k) &&\mbox{on }\Gamma_2,
    \end{aligned}\right.\end{align*}
    where $P_0:=0$. Compare with (\ref{original2}).
\end{Rem}

\medskip

In this paper, we demonstrate solvability (Proposition \ref{prop_weak})
and stability (Theorem \ref{thm_stab}) of the projection method
and establish error estimates in suitable norms 
(Theorems \ref{thm_conv1} and \ref{thm_conv2}).
Furthermore, we prove the existence of a weak solution of (\ref{original})
with a different approach than \cite{Bernard03,KC15}
by using the stability result (Corollary \ref{cor_existsol}).

The organization of this paper is as follows.
In Section 2, we introduce the notations used in this work, 
the weak formulations of the Navier--Stokes equations (\ref{original2}),
and the projection method (\ref{step1}), (\ref{step2_1}), and (\ref{step2_2}).
We also prove the existence of the weak solution to the scheme
and provide the main results.
Section \ref{sec_proofs} is devoted to proving that the solution to the scheme 
is bounded in suitable norms and converges to 
the solution to (\ref{original2}) in a strong topology as $\tau\rightarrow0$.
We also establish error estimates in suitable norms 
between the solutions to the Navier--Stokes equations
and the projection method.
In Section \ref{sec_num},
we show a numerical example of the projection method
and the numerical errors between the Navier--Stokes equations
and the projection method using the P2/P1 finite element method.
We conclude this paper with several comments on future works
in Section \ref{sec_conclusion}.
In the Appendix, four lemmas in Section \ref{sec_pre} and \ref{sec_proofs} 
and the existence result of the weak solution to (\ref{original2}) 
are proved.

\section{Preliminaries}\label{sec_pre}
In this section, we introduce the notations used in this work 
and the weak formulations of the Navier--Stokes equations (\ref{original2})
and the projection method (\ref{step1}), (\ref{step2_1}), and (\ref{step2_2}).

\subsection{Notation}\label{sec_notation}

We prepare the function spaces and the notation to be used throughout the paper.
Let $\Omega$ be a bounded Lipschitz domain in $\R^d~(d=2,3)$.
For the boundary $\Gamma=\partial\Omega$, 
we assume that there exist two relatively open subsets
$\Gamma_1,\Gamma_2$ of $\Gamma$ such that
$\Gamma_2$ has a finite number of connected components
that are piecewise $C^{1,1}$-class and
\[
  |\Gamma\setminus(\Gamma_1\cup\Gamma_2)|=0,\quad
  |\Gamma_1|,|\Gamma_2|> 0,\quad
  \Gamma_1\cap\Gamma_2=\emptyset,\quad
  \mathring{\overline{\Gamma_1}}=\Gamma_1,\quad
  \mathring{\overline{\Gamma_2}}=\Gamma_2,
\]
where $\overline{A}$ is the closure of $A\subset\Gamma$ with respect to $\Gamma$,
$\mathring{A}$ is the interior of $A$ with respect to $\Gamma$,
and $|A|$ is the $(d-1)$-dimensional Hausdorff measure of $A$. 
For an integer $m\ge 1$ and a real number $p\in[1,\infty]$,
we use the usual Lebesgue spaces $L^p(\Omega)$ and Sobolev spaces $W^{m,p}(\Omega)$
together with their standard norms
and write $H^0(\Omega):=\Lo,H^m(\Omega):=W^{m,2}(\Omega)$.
We use the same notation $(\cdot,\cdot)$ to represent 
the $\Lo$ inner product for scalar-, vector-, and matrix-valued functions.
For a normed space $X$, the dual pairing between X and the dual space $X^*$
is denoted by $\langle\cdot,\cdot\rangle_X$,
and we simply write $L^2(X)$ and $L^\infty(X)$ 
as $L^2(0,T;X)$ and $L^\infty(0,T;X)$, respectively.
$\mathcal{D}'(\Omega)$ denotes the space of distributions on $\Omega$.

We use the following notations:
\[\begin{aligned}
    \|\cdot\|_m&\mbox{: the norms in }H^m(\Omega),\\
    \Cinfb&:=\{\varphi|_\Omega~|~\varphi:\R^d\rightarrow\R\mbox{ is infinitely differentiable}\}\\
    \Cinfo&:=\{\varphi\in\Cinfb~|~\mbox{supp}(\varphi)\subset\Omega\}\\
    H&:=\{\varphi\in\Hb^d~|~\varphi=0\mbox{ on }\Gamma_1, 
        \varphi\times n=0\mbox{ on }\Gamma_2\},\\
    V&:=\{\varphi\in H~|~\dvg\varphi=0\mbox{ in }\Lo\},\\
    \Hsecond &:=\{\psi\in\Hb~|~\psi=0\mbox{ on }\Gamma_2\}.
\end{aligned}\]
The dual spaces $H^*$ and $V^*$ are equipped with the dual norms
\[
    \|F\|_{H^*}
    :=\sup_{0\ne\varphi\in H}\frac{\langle F,\varphi\rangle_H}{\|\varphi\|_1},\qquad
    \|G\|_{V^*}
    :=\sup_{0\ne\varphi\in V}\frac{\langle G,\varphi\rangle_V}{\|\varphi\|_1},
\]
respectively, for all $F\in H^*$ and $G\in V^*$.

We also use the Lebesgue space $L^2(\Gamma)$ and Sobolev space $H^{1/2}(\Gamma)$
defined on $\Gamma$.
The norm $\|\eta\|_{H^{1/2}(\Gamma)}$ is defined by
\[
    \|\eta\|_{H^{1/2}(\Gamma)}
    :=\left(\|\eta\|^2_{L^2(\Gamma)}+\int_\Gamma \int_\Gamma
    \frac{|\eta(x)-\eta(y)|^2}{|x-y|^d}ds(x)ds(y)\right)^{1/2},
\]
where $ds$ denotes the surface measure of $\Gamma$.
For function spaces defined on $\Gamma_2$, 
we use $L^2(\Gamma_2)$ and $H^{1/2}(\Gamma_2)$.

We define a bilinear form $a:\Hb^d\times\Hb^d\rightarrow\R$, 
a seminorm on $\Hb^d$, for $u,v\in \Hb^d$,
\[\begin{aligned}
    a(u,v):=\into(\dvg u)(\dvg v)dx+\into(\rot u)\cdot(\rot v)dx,\qquad
    \|u\|_a:=\sqrt{a(u,u)}.
\end{aligned}\]
Let $p_d$ be 
\[
    p_d
    :=\left\{\begin{array}{ll}
        2+\eps &\mbox{if }d=2,\\
        3 &\mbox{if }d=3,
    \end{array}\right. 
\]
where $\eps>0$ is arbitrarily small.
It follows from the Sobolev embeddings that 
$\Hb \subset L^{p_d}(\Omega)$ and the embedding is continuous
\cite[Theorem III.2.34]{BF13}.
We define a trilinear operator
$d:L^{p_d}(\Omega)^d\times H\times H\rightarrow\R$
for $u\in L^{p_d}(\Omega)^d$ and $v,w\in H$,
\[
    d(u,v,w)
    :=\into u\cdot(\nabla\times(v\times w))dx.
\]
We note that for all $u\in\Hb^d$ and $v,w\in H$,
\[
    d(u,v,w)
    =-\int_{\Gamma}(u\times\nu)\cdot(v\times w)ds
    +\into((\nabla\times u)\times v)\cdot wdx
    =\into D^1(u,v)\cdot wdx.
\]

For two sequences $(x_k)_{k=0}^N$ and $(y_k)_{k=1}^N$ 
in a Banach space $E$, we define a piecewise linear interpolant 
$\hat{x}_\tau\in W^{1,\infty}(0,T;E)$ of $(x_k)_{k=0}^N$ 
and a piecewise constant interpolant 
$\bar{y}_\tau\in L^\infty(0,T;E)$ of $(y_k)_{k=1}^N$, respectively, by
\[\begin{aligned}
    \hat{x}_\tau(t)&:=x_{k-1}+\frac{t-t_{k-1}}{\tau}(x_k-x_{k-1})
    &&\text{for }t\in[t_{k-1},t_k] \mbox{ and }k=1,2,\ldots,N,\\
    \bar{y}_\tau(t)&:=y_k
    &&\text{for }t\in(t_{k-1},t_k] \mbox{ and }k=1,2,\ldots,N.
\end{aligned}\]
We define a backward difference operator by 
\[
    D_\tau x_k:=\frac{x_k-x_{k-1}}{\tau},\qquad
    D_\tau y_l:=\frac{y_l-y_{l-1}}{\tau}
\]
for $k=1,2,\ldots,N$ and $l=2,3,\ldots,N$.
Then, the sequence $(D_\tau x)_k:=D_\tau x_k$ satisfies 
$\frac{\partial\hat{x}_\tau}{\partial t}=(\overline{D_\tau x})_\tau$
on $(t_{k-1},t_k)$ for all $k=1,2,\ldots,N$.
For a function $F\in C([0,T];E)$, we define $F_\tau\in L^\infty(0,T;E)$
as the piecewise constant interpolant of $(F(t_k))_{k=1}^N$, i.e.,
\[
    F_\tau(t):=F(t_k)
    \qquad\mbox{for }t\in(t_{k-1},t_k] \mbox{ and }k=1,2,\ldots,N.
\]

\subsection{Preliminary results}
Let $\gamma_0\in B(\Hb, H^{1/2}(\Gamma))$ be
the standard trace operator.
It is known that (see e.g. \cite[Theorem 1.2]{Temam})
there exists a linear continuous operator
$\gamma_n :M:=\{\varphi\in\Lo^d~|~\dvg\varphi\in\Lo\}
\rightarrow H^{-1/2}(\Gamma)$ such that
$\gamma_n u=u\cdot n |_\Gamma$ for all $u\in C^\infty(\overline{\Omega})^d$,
where $H^{-1/2}(\Gamma):=H^{1/2}(\Gamma)^*$.
Then, the following generalized Gauss divergence formula holds:
\[
    \into v\cdot\nabla q dx + \into(\dvg v)q dx
    = \langle\gamma_n v,\gamma_0 q\rangle_{H^{1/2}(\Gamma)}
    \qquad\text{ for all }v\in M,q\in\Hb.
\]
The composition of the trace operator $\gamma_0$ and 
the restriction $H^{1/2}(\Gamma)\rightarrow H^{1/2}(\Gamma_2)$ 
is denoted by $\psi\mapsto\psi|_{\Gamma_2}$. 
This map is continuous from $\Hb$ to $H^{1/2}(\Gamma_2)$.
The kernel of this map is $H^1_{\Gamma_2}(\Omega)$.
We simply write $\psi$ instead of $\psi|_{\Gamma_2}$
when there is no ambiguity.

We recall the following lemmas that are necessary for
the existence and the uniqueness of a solution to the Stokes problem.

\begin{Lem}\label{Stokes}
  {\rm \cite[Corollary 4.1]{Girault}}
Let $(X,\|\cdot\|_X)$ and $(Q,\|\cdot\|_Q)$ be two real Hilbert spaces.
Let $a_X:X\times X\rightarrow\R$ and $b:X\times Q\rightarrow\R$ be
bilinear and continuous maps and let $F\in X^*$.
If there exist two constants $\alpha>0$ and $\beta>0$ such that
for all $v\in V$ and $q\in Q$,
\[
    a_X(v,v)\ge\alpha\|v\|^2_X,\qquad
    \sup_{0\ne v\in X}\frac{b(v,q)}{\|v\|_{X}}\ge \beta\|q\|_Q
\]
where $V=\{v\in X~|~b(v,q)=0\mbox{ for all }q\in Q\}$,
then there exists a unique solution $(u,p)\in X\times Q$ to the following problem:
\begin{align*}\left\{\begin{aligned}
    &a_X(u,v)+b(v,p)=F(v)&\mbox{for all }v\in X,\\
    &b(u,q)=0&\mbox{for all }q\in Q.
\end{aligned}\right.\end{align*}
Furthermore, there exists a constant $c>0$ independent of $F$ such that 
\[
    \|u\|_X+\|p\|_Q\le c\|F\|_{X^*}.
\]
\end{Lem}

\begin{Lem}\label{lem_infsup}
    {\rm \cite[proof of Theorem 2.1]{BCM91}}
    There exists a constant $c=c(\Omega,\Gamma_1,\Gamma_2)$ $>0$ 
    such that for all $q\in\Lo$,
    \[
        \|q\|_0
        \le c\sup_{0\ne\varphi\in H}\frac{|(q,\dvg\varphi)|}{\|\varphi\|_1}.
    \]
\end{Lem}

The following embedding theorem show
the continuity and the coercivity of the bilinear form 
$a:H\times H\rightarrow\R$.

\begin{Lem}\label{lem_a}
    There exists a constant $c_a=c_a(\Omega,\Gamma_1,\Gamma_2)>0$ 
    such that for all $\varphi_1,\varphi_2,\varphi\in H$,
    \[\begin{aligned}
        a(\varphi_1,\varphi_2)\le\|\varphi_1\|_a\|\varphi_2\|_a
        \le c_a\|\varphi_1\|_1\|\varphi_2\|_1,\qquad
        \frac{1}{c_a}\| \varphi\|_1^2\le \|\varphi\|_a^2.
    \end{aligned}\] 
\end{Lem}
The first inequality holds from the Cauchy--Schwarz inequality.
For the proof of the second inequality,
see \cite[Lemma 2.11]{BRY15} and \cite[Lemma 5]{KN99}.
The following embedding theorem is called the Poincar\'e inequality.

\begin{Lem}[Poincar\'e's inequality]\label{Poincare}
  {\rm \cite[Lemma 3.1]{Girault}}
  There exists a constant $c=c(\Omega,\Gamma_1,$ $\Gamma_2)>0$ 
  such that for all $\varphi\in\Hsecond$,
  \[
    \|\varphi\|_0\le c\|\nabla\varphi\|_0.
  \] 
\end{Lem}

We prepare the following lemma to use the Aubin--Nitsche trick.

\begin{Lem}\label{lem_stokes}
    We define an operator 
    $T:\Lo^d\ni e\mapsto(w,r)\in H\times\Lo$ as follows:
    \begin{align}\label{eq_stokes}\left\{\begin{array}{ll}
        {\DS a(w,\varphi)-(r,\dvg\varphi)=(e,\varphi)}
        &\mbox{for all }\varphi\in H,\\
        \dvg w=0 &\mbox{in }\Lo.
    \end{array}\right.\end{align}
    Then, $T$ is a linear and continuous operator and
    there exists a constant $c=c(\Omega,\Gamma_1,\Gamma_2)>0$ 
    such that for all $e\in\Lo^d$ and $(w,r)=T(e)$,
    \[
        \|w\|_1+\|r\|_0\le c\|e\|_{H^*},\qquad
        \frac{1}{c}\|e\|_{V^*}
        \le \|w\|_1 
        \le c\|e\|_{V^*}.
    \]
\end{Lem}
By Lemmas \ref{Stokes}, \ref{lem_infsup}, and \ref{lem_a},
the operator $T$ is well-posed and continuous.
See the Appendix for the proof of the inequalities. 
Next, we show the following two lemmas for the operator $d$.

\begin{Lem}\label{lem_outflow}
    It holds that for all 
    $u\in L^{p_d}(\Omega)^d,v,v_1,v_2\in\Hb^d$
    \[
        d(u,v,v)=0,\qquad
        d(u,v_1,v_2)=-d(u,v_2,v_1).
    \]
\end{Lem}

By the definition of the operator $d$, 
it is easy to check Lemma \ref{lem_outflow}.
\begin{Lem}\label{lem_nonlin}
  There exists a constant $c_d=c_d(\Omega,\Gamma_1,\Gamma_2)>0$
  such that
  \[
      d(u,v,w)\le \left\{\begin{array}{l}
          c_d\|u\|_{L^{p_d}}\|v\|_1\|w\|_1
          \text{ for all }u \in L^{p_d}(\Omega)^d, v, w \in H,\\
          c_d\|u\|_0\|v\|_1\|w\|_2
          \text{ for all }u \in L^{p_d}(\Omega)^d, v \in H, 
          w \in H \cap H^2(\Omega)^d,\\
          c_d\|u\|_1\|v\|_1\|w\|_1
          \text{ for all }u \in \Hb^d, v, w \in H,\\
          c_d\|u\|_1\|v\|_2\|w\|_0
          \text{ for all }u \in \Hb^d, v \in H \cap H^2(\Omega)^d, 
          w \in H\\
          c_d\|u\|_2\|v\|_1\|w\|_0
          \text{ for all }u \in H^2(\Omega)^d, v, w \in H.
      \end{array}\right.
  \]
\end{Lem}

See the Appendix for the proof.
Finally, we recall the discrete Gronwall inequality.

\begin{Lem}\label{lem_gronwall}
    {\rm \cite[Lemma 5.1]{HR90}}
  Let $\tau,\beta>0$ and let nonnegative sequences 
  $(a_k)^N_{k=0}$, $(b_k)^N_{k=0}$, $(c_k)^N_{k=0}$, 
  $(\alpha_k)^N_{k=0} \subset \{x\in\R ~|~ x \ge 0\}$ satisfy that
  \[
      a_n+\tau\sum^m_{k=0}b_k
      \le \tau\sum^m_{k=0}\alpha_k a_k+\tau\sum^m_{k=0}c_k+\beta
      \qquad\mbox{for all }m=0,1,\ldots,N.
  \]
  If $\tau\alpha_k<1$ for all $k=0,1,\ldots,N$, then we have
  \[
      a_n+\tau\sum^m_{k=0}b_k
      \le e^C\left(\tau\sum^m_{k=0}c_k+\beta\right)
      \qquad\mbox{for all }m=0,1,\ldots,N,
  \]
  where $C:=\tau\sum^N_{k=0}\frac{\alpha_k}{1-\tau\alpha_k}$.
\end{Lem}

\subsection{Weak formulations of (\ref{original2}), 
(\ref{step1}), (\ref{step2_1}), and (\ref{step2_2})}\label{sec_weak}

We assume $\nu=\rho=1$ and the following conditions for $f,p^b,$ and $u_0$:
\begin{align}\label{cond}
    f\in L^2(H^*),\quad 
    p^b\in L^2(\Hb),\quad
    u_0\in L^{p_d}(\Omega)^d.
\end{align}
To define weak formulations of the Navier--Stokes equations
(\ref{original2}) and the projection method
(\ref{step1}), (\ref{step2_1}), and (\ref{step2_2}), 
we prepare the following equation:

\begin{Prop}\label{prop_weak}
    It holds that for all $u\in H^2(\Omega)$ and $\varphi\in H$, 
    \begin{align}\label{eq_asmooth}
        -(\Delta u,\varphi)
        =a(u,\varphi)-\int_{\Gamma_2}(\dvg u)\varphi\cdot n ds.
    \end{align}
\end{Prop}
\noindent{\it Proof.}
It holds that
$
    -\Delta u=\rot(\rot u)-\nabla(\dvg u)
$
for all $u\in C^2(\overline{\Omega})^d$.
Hence, we have for all $u\in C^2(\overline{\Omega})^d$
and $\varphi\in C^1(\overline{\Omega})^d$, 
\[
    (-\Delta u,\varphi)
    =a(u,\varphi)+\int_\Gamma (\rot u)\cdot(\varphi\times n)ds
    -\int_\Gamma (\dvg u)\varphi\cdot n ds,
\]
which also holds for all $\varphi\in H^2(\Omega)$ and $\psi\in\Hb$
since the two spaces $C^2(\overline{\Omega})$ and $C^1(\overline{\Omega})$
are dense in $H^2(\Omega)$ and $\Hb$, respectively.
By the definition of $H$,
equation (\ref{eq_asmooth}) holds for all $u\in H^2(\Omega)$ and $\varphi\in H$.
\qed\medskip

By Proposition \ref{prop_weak} and the Gauss divergence formula,
it holds that for all $u\in H^2(\Omega)^d,P\in\Hb$, and $\varphi\in V$
with $\dvg u=0$ in $\Hb$,
\[
    (D(u,u)-\Delta u+\nabla P,\varphi)
    =a(u,\varphi)+d(u,u,\varphi)-(P,\dvg\varphi)
    +\int_{\Gamma_2}P\varphi\cdot nds.
\]
Hence, a weak formulation of (\ref{original2}) is as follows: 
Find $u\in L^2(\Hb^d)$ and $P\in L^1(\Lo)$ such that 
$\frac{\partial u}{\partial t}\in L^1(H^*)$
and for all $\varphi\in H$,
\begin{align}\tag{NS}\label{weak_original}\left\{\begin{array}{l}
    {\DS \left\langle\frac{\partial u}{\partial t},\varphi\right\rangle_H
    + a(u,\varphi) + d(u,u,\varphi) - (P,\dvg\varphi)
    = \langle f,\varphi\rangle_H
    - \int_{\Gamma_2}p^b\varphi\cdot n ds},\\[8pt]
    \dvg u=0\qquad\mbox{in }\Lo,
\end{array}\right.\end{align}
in $L^1(0,T)$, with $u(0)=u_0$ on $H^*$. In main convergence theorems 
(Theorems \ref{thm_conv1} and \ref{thm_conv2}), 
we assume that (\ref{weak_original}) has a unique solution 
and that the solution is as smooth as needed. 

On the other hand, by Proposition \ref{prop_weak},
we have for all $u\in H^2(\Omega)^d$ and $\varphi\in H$,
\[
    (D(u,u)-\Delta u,\varphi)
    =a(u,\varphi)+d(u,u,\varphi)
    -\int_{\Gamma_2}(\dvg u)\varphi\cdot nds.
\]
Hence, a weak formulation of (\ref{step1}), (\ref{step2_1}), and (\ref{step2_2})
with the initial datum $u_0(=:\upre_0)$ is as follows: 
\begin{Prob}\label{def_FS}
    Let $(f_k)_{k=1}^N \subset H^*$ and $(p^b_k)_{k=1}^N \subset \Hb$.
    For all $k=1,2,\ldots,$ $N$, find $(\upre_k,P_k,$\\ $u_k)\in H\times\Hb\times\Lo^d$
    such that $P_k-p^b_k\in\Hsecond$ and 
    for all $\varphi\in H$ and $\psi\in\Hsecond$,
    \begin{align}\tag{PM}\label{weak_FS}\left\{\begin{array}{l}
        {\DS \frac{1}{\tau}(\upre_k-u_{k-1},\varphi)
        + a(\upre_k,\varphi) + d(\upre_{k-1},\upre_k,\varphi)
        = \langle f_k,\varphi\rangle_H}\\[8pt]
        {\DS \tau(\nabla P_k,\nabla\psi)
        = -(\dvg\upre_k,\psi)}\\[4pt]
        u_k=\upre_k-\tau\nabla P_k \text{ in }\Lo^d.
    \end{array}\right.\end{align}
\end{Prob}

\begin{Rem}\label{rem_fpb}
    For $f\in L^2(H^*)$ and $p^b\in L^2(\Hb)$, 
    we set for all $k=1,2,\ldots,N$,
    \begin{align}\label{def_fpb}
        f_k := \frac{1}{\tau}\int_{t_{k-1}}^{t_k} f(s)ds,\qquad 
        p^b_k := \frac{1}{\tau}\int_{t_{k-1}}^{t_k} p^b(s)ds.
    \end{align} 
    Here, it holds that $\bar{f}_\tau \in L^2(H^*)$
    and $\bar{p}_\tau^b \in L^2(\Hb)$:
    \[
        \|\bar{f}_\tau\|_{L^2(H^*)}
        \le \|f\|_{L^2(H^*)},\qquad
        \|\bar{p}^b_\tau\|_{L^2(H^1)}
        \le \|p^b\|_{L^2(H^1)}.
    \]
    In Theorems \ref{thm_conv1} and \ref{thm_conv2},
    we assume $f \in C([0,T];H^*), p^b \in C([0,T];\Hb)$
    to use $f(t_k)$ and $p^b(t_k)$ for all $k=1,2,\ldots, N$
    (Hypothesis \ref{hyp_reg_sol}).
    Then, we set for all $k=1,2,\ldots,N$,
    \begin{align*}
        f_k := f(t_k),\qquad 
        p^b_k := p^b(t_k),
    \end{align*} 
    which implies that $\bar{f}_\tau = f_\tau \in L^2(H^*)$
    and $\bar{p}^b_\tau = p^b_\tau \in L^2(\Hb)$.
\end{Rem}

\medskip
We show the existence and uniqueness of the solution to (\ref{weak_FS})
in the following proposition.

\begin{Prop}
    For all $(f_k)_{k=1}^N \subset H^*$, 
    $(p^b_k)_{k=1}^N \subset \Hb^d$, and $u_0\in L^{p_d}(\Omega)^d$,
    Problem \ref{def_FS} has a unique solution.
\end{Prop}

\noindent{\it Proof.}
By Lemmas \ref{lem_a}, \ref{lem_outflow}, and \ref{lem_nonlin},
if $\upre_{k-1} \in L^{p_d}(\Omega)^d$ are known, 
then it holds that for all $v,\varphi\in H$,
\[\begin{aligned}
    \frac{1}{\tau}(v,\varphi) + a(v,\varphi) + d(\upre_{k-1},v,\varphi)
    &\le \left(\frac{1}{\tau} + c_a 
    + c_d\|\upre_{k-1}\|_{L^{p_d}}\right)\|v\|_1\|\varphi\|_1,\\
    \frac{1}{\tau}(v,v) + a(v,v) + d(\upre_{k-1},v,v)
    &\ge \frac{1}{c_a}\|v\|_1^2,
\end{aligned}\]
which implies that the mapping $H\times H\ni(v,\varphi)\mapsto 
\frac{1}{\tau}(v,\varphi) + a(v,\varphi) + d(\upre_{k-1},$ $v,\varphi)\in\R$
is a continuous and coercive bilinear form.
On the other hand, if $u_{k-1} \in \Lo^d$, then the mapping 
$H\ni\varphi\mapsto
\langle f(t_k),\varphi\rangle_H + \tau^{-1}(u_{k-1},\varphi)\in\R$
is a functional on $H$.
By the Lax--Milgram theorem, there exists a unique solution 
$\upre_k \in H \subset L^{p_d}(\Omega)^d$ 
to the first equation of (\ref{weak_FS}).
Since $\dvg\upre_k \in \Lo$, by the Poincar\'e inequality and 
the Lax--Milgram theorem, the second equation of (\ref{weak_FS}) 
also has a unique solution $P_k \in \Hb$.
Furthermore, we obtain 
$u_k := \upre_k - \tau\nabla P_k \in \Lo^d$.
Therefore, since $u_0(=\upre_0) \in L^{p_d}(\Omega)^d$, 
(\ref{weak_FS}) has a unique solution 
$(\upre_k,P_k,u_k)_{k=1}^N \subset H\times\Hb\times\Lo^d$.
\qed

\begin{Rem}\label{rem_ortho}
    The function space $\Lo^d$ has the following 
    orthogonal decomposition:
    \[
        \Lo^d = U \oplus \nabla(\Hsecond),
    \]
    where $U := \{\varphi\in\Lo^d ~|~ \dvg\varphi=0 \text{ in }\Lo,
        \langle \gamma_n\varphi,\psi\rangle_{H^{1/2}(\Gamma)} = 0 
        \text{ for all } \psi\in\Hsecond\}$ \cite[Proposition 4.1]{GQ98}.
    By the second and third equation of (\ref{weak_FS})
    and the Gauss divergence formula,
    it holds that for all $k=1,2,\ldots,N$ and $\psi\in\Hsecond$,
    \[
        (u_k,\nabla\psi)
        = (\upre_k,\nabla\psi) - \tau(\nabla P_k,\nabla\psi)
        = - (\dvg\upre_k,\psi) - \tau(\nabla P_k,\nabla\psi)
        = 0,
    \]
    which implies that $u_k \in U$.
    Since the third equation of (\ref{weak_FS}) is equivalent to 
    \[
        \upre_k - \tau\nabla p^b(t_k)
        = u_k + \tau\nabla (P_k-p^b(t_k)) 
        \qquad\text{ in }\Lo^d,
    \]
    Step 2 ((\ref{step2_1}) and (\ref{step2_2})) 
    is the projection of $\upre_k - \tau\nabla p^b(t_k)$ 
    to the divergence-free space $U$.
\end{Rem}

\begin{Rem}\label{rem_equiv}
    By replacing $u_{k-1}$ in the first equation of (\ref{weak_FS})
    with the third equation of (\ref{weak_FS}) at the previous step,
    it holds that for all $k=1,2,\ldots,N$,
    $\varphi\in H$, and $\psi\in\Hsecond$,
    \begin{align*}\left\{\begin{array}{l}
        {\DS \frac{1}{\tau}(\upre_k-\upre_{k-1},\varphi)
        + a(\upre_k,\varphi) + d(\upre_{k-1},\upre_k,\varphi)
        + (\nabla P_{k-1},\varphi)
        = \langle f_k,\varphi\rangle_H}\\[8pt]
        {\DS \tau(\nabla P_k,\nabla\psi)
        = -(\dvg\upre_k,\psi)}
    \end{array}\right.\end{align*}
    where $P_0:=0$ (cf. \cite{Rannacher92}). 
     Ones can calculate $(\upre_k,P_k)_{k=1}^N$
    without the velocity $(u_k)_{k=1}^M$. 
    Since the calculation $u_k=\upre_k-\tau\nabla P_k$ is not used,
    this formulation is suitable for 
    numerical calculations such as the finite element method 
    (see Section \ref{sec_num}).

    On the other hand, 
    by replacing $\upre_k$ in the first term of the first equation of (\ref{weak_FS})
    with the third equation of (\ref{weak_FS}) at the same step,
    it holds that for all $k=1,2,\ldots,N$,
    $\varphi\in H$, and $\psi\in\Hsecond$,
    \begin{align}\label{weak_FS_equiv2}\left\{\begin{array}{l}
        {\DS \frac{1}{\tau}(u_k-u_{k-1},\varphi)
        + a(\upre_k,\varphi) + (\nabla P_k,\varphi)
        = \langle f_k,\varphi\rangle_H - d(\upre_{k-1},\upre_k,\varphi)}\\[8pt]
        {\DS \tau(\nabla P_k,\nabla\psi) + (\dvg\upre_k,\psi)
        = 0}\\[4pt]
        u_k=\upre_k-\tau\nabla P_k \text{ in }\Lo^d.
    \end{array}\right.\end{align}
    This formulation is helpful to prove stability and convergence results.
\end{Rem}

\subsection{Main theorems for stability and convergence}

We show the stability of the projection method (\ref{weak_FS})
and establish error estimates in suitable norms 
between the solutions to the Navier--Stokes equations 
(\ref{weak_original}) and the projection method (\ref{weak_FS}).

\begin{Th}\label{thm_stab}
    Under the condition (\ref{cond}),
    we set $f_k \in H^*$ and $p^b_k \in \Hb^d$ as (\ref{def_fpb})
    for all $k=1,2,\ldots,N$. 
    Then, there exists a constant $c>0$ independent of $\tau$ such that
    \[\begin{aligned}
        &\quad\, \|\bar{u}_\tau\|_{L^\infty(L^2)}
        + \|\bar{u}^*_\tau\|_{L^\infty(L^2)}
        + \|\bar{u}^*_\tau\|_{L^2(H^1)}
        + \frac{1}{\sqrt{\tau}}\|\bar{u}_\tau-\bar{u}^*_\tau\|_{L^2(L^2)}\\
        &\le c\left(\|u_0\|_0 + \|f\|_{L^2(H^*)} 
        + \|p^b\|_{L^2(H^1)}\right).
    \end{aligned}\]
\end{Th}

For a convergence theorem, we assume:

\begin{Hyp}\label{hyp_reg_sol}
    The solution $(u,P)$ to (\ref{weak_original}) satisfies
    \begin{align*}
        u\in C([0,T];H\cap H^2(\Omega)^d)\cap H^1(\Lo^d)\cap H^2(H^*),\qquad
        P\in C([0,T];\Hb).
    \end{align*}
    We also assume $f \in C([0,T];H^*)$ and $p^b \in C([0,T];\Hb)$
    and set in Problem \ref{def_FS} for all $k=1,2,\ldots,N$,
    \begin{align*}
        f_k := f(t_k),\qquad 
        p^b_k := p^b(t_k).
    \end{align*} 
\end{Hyp}

\begin{Th}\label{thm_conv1}
    Under Hypothesis \ref{hyp_reg_sol},
    there exist two constants $c,\tau_0>0$ independent of $\tau$ 
    such that for all $0<\tau<\tau_0$,
    \[\begin{aligned}
        \|u-\bar{u}_\tau\|_{L^\infty(L^2)}
        +\|u-\bar{u}^*_\tau\|_{L^\infty(L^2)}
        +\|u-\bar{u}^*_\tau\|_{L^2(H^1)}
        &\le c\sqrt{\tau},\\
        \|\bar{u}_\tau-\bar{u}^*_\tau\|_{L^2(L^2)}
        &\le c\tau.
    \end{aligned}\]
\end{Th}

\begin{Rem}
    For regularity of the solution $(u,P)$ to (\ref{weak_original}),
    see \cite[Theorem 1.3]{Bernard03} 
    and \cite[Theorems 4.2 and 4.3]{Kim15}.
    In the case of the homogeneous Dirichlet boundary condition
    on the whole boundary $\Gamma$, high regularity properties of 
    the solution to the Navier--Stokes equations are proved 
    in \cite[Theorem V.2.10]{BF13}.
\end{Rem}

\begin{Rem}\label{rem_smooth}
    If $u\in C([0,T];H\cap H^2(\Omega)^d)$,
    then $|u|^2\in C([0,T];\Hb)$, and hence,
    $p\in C([0,T];\Hb)$ is equivalent to
    $P=p+\frac{1}{2}|u|^2\in C([0,T];\Hb)$.
\end{Rem}

%

Furthermore, we assume the following regularity assumptions:

\begin{Hyp}[Regularity of the Stokes problem]\label{hyp_stokes}
    There exists a constant $c = c(\Omega,\Gamma_1,$\\ $\Gamma_2)>0$ such that
    \[
        \|w\|_2 + \|r\|_1 \le c\|e\|_0.
    \]
    for all $e \in \Lo^d$ and $(w, r) = T(e)$.
\end{Hyp}

\begin{Hyp}\label{hyp_reg_sol2}
    The solution $(u,P)$ to (\ref{weak_original}) satisfies
    \[
        u \in H^1(\Hb^d) \cap H^2(\Lo^d) \cap H^3(H^*),\qquad
        P \in H^1(\Hb).
    \]
\end{Hyp}

Then we can improve the convergence rate:

\begin{Th}\label{thm_conv2}
    Under Hypothesis \ref{hyp_reg_sol} and \ref{hyp_stokes},
    there exist two constants $\tau_1,c>0$ independent of $\tau$ 
    such that for all $0<\tau<\tau_1$,
    \[
        \|u-\bar{u}_\tau\|_{L^2(L^2)}
        + \|u-\bar{u}^*_\tau\|_{L^2(L^2)}
        \le c\tau.
    \]
    Furthermore, if we also assume Hypothesis \ref{hyp_reg_sol2},
    then there exist two constants $\tau_2,c>0$ independent of $\tau$ 
    such that for all $0<\tau<\tau_2(\le \tau_1)$,
    \[
        \|P-\bar{P}_\tau\|_{L^2(L^2)}
        \le c\sqrt{\tau}.
    \]
\end{Th}

\begin{Rem}
    Hypothesis \ref{hyp_stokes} holds, e.g., 
    if $\Omega$ is of class $C^{2,1}$ \cite[Theorem 1.2]{Bernard02}.
\end{Rem}

\subsection{Main result for existence of a weak solution to (\ref{original2})}

Using Theorem \ref{thm_stab}, we prove that there exists a solution 
to a weak formulation of (\ref{original2}) weaker than (\ref{weak_original}).
Putting $\varphi:=v\in V$ in the first equation of (\ref{weak_original}), 
we obtain the following equation: for all $v\in V$,
\begin{align}\label{weak_original2}
    \left\langle\frac{\partial u}{\partial t},v\right\rangle_V
    + a(u,v) + d(u,u,v)
    =\langle f,v\rangle_H - \int_{\Gamma_2}p^b v\cdot n ds
\end{align}
in $L^1(0,T)$.

\begin{Cor}\label{cor_existsol}
    Under the condition (\ref{cond}),
    there exists a solution $u\in L^2(V)\cap L^\infty(\Lo^d)\cap C([0,T];V^*)$
    to (\ref{weak_original2}) with $u(0)=u_0$
    such that $\frac{\partial u}{\partial t}\in L^{4/p_d}(V^*)$.
\end{Cor}

\begin{Rem}
    For $f\in L^2(\Lo^d)$, local existence and uniqueness of 
    a weak solution $u$ to (\ref{weak_original2}) with $u_0\in H$ 
    are proved in \cite[Theorem 1.3]{Bernard03}.
    By \cite[Lemma 4]{KN99}:
    \[
        a(u,v)
        = \sum^d_{i,j=1}\into\frac{\partial u_i}{\partial x_j}\frac{\partial v_i}{\partial x_j} dx
        + \int_{\Gamma_2}\kappa u\cdot vds
        \quad\text{for all }u,v\in H,
    \]
    where $\kappa:=\dvg n=(d-1)\times $(mean curvature)
    (cf. Remark \ref{rem_KN}), (\ref{weak_original2}) is equivalent to 
    \begin{align}\label{weak_CK}\begin{aligned}
        &\left\langle\frac{\partial u}{\partial t},v\right\rangle_V
        +\sum^d_{i,j=1}\into\frac{\partial u_i}{\partial x_j}\frac{\partial v_i}{\partial x_j} dx+\int_{\Gamma_2}\kappa u\cdot vds+d(u,u,v)\\
        =&\langle f,v\rangle_H-\int_{\Gamma_2}p^b v\cdot n ds
        \quad\text{for all }v\in V
    \end{aligned}\end{align}
    in $L^1(0,T)$.
    It is known \cite[Theorem 5.1]{KC15} that 
    there exists a weak solution 
    $u$ 
    to (\ref{weak_CK}) with $u_0\in U$, 
    where $U$ is defined in Remark \ref{rem_ortho}.
\end{Rem}

\section{Proofs}\label{sec_proofs}
In this section, we prove that the solution to (\ref{weak_FS})
is bounded in suitable norms (Theorem \ref{thm_stab}) 
and error estimates (Theorems \ref{thm_conv1} and \ref{thm_conv2})
in suitable norms between the solutions 
to (\ref{weak_original}) and (\ref{weak_FS}).

\subsection{Stability}\label{sec_stab}

We prepare the following useful lemma for the proofs of 
Theorems \ref{thm_stab}, \ref{thm_conv1}, and \ref{thm_conv2}.

\begin{Lem}\label{lem_gen}
    Let $v_0 \in \Lo^d, 
    (F_k, G_k, Q_k)^N_{k=1} \subset H^* \times H^* \times \Hb$
    and let $(\vpre_k, v_k, q_k)_{k=1}^N \in H \times \Lo^d \times \Hb$
    satisfy that for all $k = 1,2,\ldots,N$, 
    $\varphi \in H$, and $\psi \in \Hsecond$,
    \begin{align}\label{eq_gen}\left\{\begin{aligned}
        &\frac{1}{\tau}(v_k - v_{k-1}, \varphi)
        + a(\vpre_k, \varphi) - (q_k, \dvg\varphi)
        = \langle F_k + G_k,\varphi\rangle_H,\\
        &\tau(\nabla q_k, \nabla\psi) + (\dvg\vpre_k, \psi)
        = -\tau(\nabla Q_k,\nabla\psi),\\
        &v_k = \vpre_k - \tau\nabla (q_k+Q_k) \mbox{ in }\Lo^d.
    \end{aligned}\right.\end{align}
    If we assume that for all $\delta>0$ 
    there exist a constant $A_\delta>0$ independent of $k$ and $\tau$, 
    and a sequence $(\beta_k)^N_{k=1}\subset\R$ such that 
    \begin{align}\label{asmp_G}
        \langle G_k, \vpre_k \rangle_H
        \le \delta\|\vpre_k\|_1^2 
        + A_\delta(\|\vpre_{k-1}\|_0^2 + \beta_k^2)
        \qquad\mbox{for all }k=1,2,\ldots,N,
    \end{align}
    where $\vpre_0:=v_0$,
    then there exist two constants $\tau_0,c>0$ independent of $\tau$ 
    such that for all $0<\tau<\tau_0$,
    \begin{align}\label{ineq_gen}\begin{aligned}
        &\quad \|\bar{v}_\tau\|_{L^\infty(L^2)}^2
        + \|\bar{v}^*_\tau\|_{L^\infty(L^2)}^2
        + \|\bar{v}^*_\tau\|_{L^2(H^1)}^2
        + \tau\left\|\frac{\partial \hat{v}_\tau}{\partial t}\right\|_{L^2(L^2)}^2
        \!\!+ \frac{1}{\tau}\|\bar{v}_\tau-\bar{v}^*_\tau\|_{L^2(L^2)}^2\\
        &\le c\left(\|v_0\|_0^2 + \|\bar{F}_\tau\|_{L^2(H^*)}^2
        + \tau\|\bar{Q}^b_\tau\|_{L^2(H^1)}^2
        + \|\bar{\beta}_\tau\|_{L^2(0,T)}^2\right).
    \end{aligned}\end{align}
    In particular, if $\langle G_k,\vpre_k\rangle_H \le 0$ 
    for all $k=1,2,\ldots,N$, then $\tau_0=T$.
\end{Lem}
\paragraph{Proof.}
Putting $\varphi:=\vpre_k$ and $\psi:=q_k$ 
and adding the two equations, we obtain for all $k=1,2,\ldots,N$,
\begin{align*}
    &\quad\, \frac{1}{\tau}(v_k-v_{k-1},\vpre_k)
    + \|\vpre_k\|_a^2 + \tau\|\nabla q_k\|_0^2
    + \tau(\nabla Q_k, \nabla q_k)\\
    &=\langle F_k + G_k, \vpre_k \rangle_H 
    \le \frac{c_a}{2}\|F_k\|_{H^*}^2 + \frac{1}{2c_a}\|\vpre_k\|_1 
    + \langle G_k, \vpre_k \rangle_H.
\end{align*}
Here, by Lemma \ref{lem_a} and the third equation of (\ref{eq_gen}), 
it holds that 
\begin{align*}
    &\quad\, \frac{1}{\tau}(v_k-v_{k-1},\vpre_k)
    + \|\vpre_k\|_a^2 + \tau\|\nabla q_k\|_0^2
    + \tau(\nabla Q_k, \nabla q_k)\\
    &= \frac{1}{\tau}(v_k-v_{k-1},v_k)
    + \frac{1}{\tau}(v_k-v_{k-1},\vpre_k-v_k)
    + \|\vpre_k\|_a^2 + \tau\|\nabla(q_k+Q_k)\|_0^2
    - \tau(\nabla Q_k, \nabla (q_k+Q_k))\\
    &\ge \frac{1}{2\tau}(\|v_k\|_0^2 - \|v_{k-1}\|_0^2 + \|v_k-v_{k-1}\|_0^2)
    - \frac{3}{8\tau}\|v_k-v_{k-1}\|_0^2
    - \frac{2}{3\tau}\|\vpre_k-v_k\|_0^2
    + \frac{1}{c_a}\|\vpre_k\|_1^2\\
    &\quad 
    + \tau\|\nabla(q_k+Q_k)\|_0^2
    - 3\tau\|\nabla Q_k\|_0^2 - \frac{\tau}{12}\|\nabla (q_k+Q_k)\|_0^2\\
    &=\! \frac{1}{2\tau}\!\left(\!\|v_k\|_0^2 - \|v_{k-1}\|_0^2
    +\! \frac{\tau^2}{4}\|D_\tau v_k\|_0^2
    +\! \frac{1}{2}\|\vpre_k-v_k\|_0^2 \!\right)\!
    +\! \frac{1}{c_a}\|\vpre_k\|_1^2\!
    -\! 3\tau\|\nabla Q_k\|_0^2. 
\end{align*}
Hence, we have for all $k=1,2,\ldots,N$,
\begin{align}\label{ineq_genp1}\begin{aligned}
    &\quad\, \|v_k\|_0^2 - \|v_{k-1}\|_0^2
    + \frac{\tau^2}{4}\|D_\tau v_k\|_0^2
    + \frac{1}{2}\|\vpre_k-v_k\|_0^2
    + \frac{\tau}{c_a}\|\vpre_k\|_1^2\\
    &\le c_a \tau\|F_k\|_{H^*}^2
    + 6\tau^2\|\nabla Q_k\|_0^2
    + 2\tau\langle G_k,\vpre_k\rangle_H.
\end{aligned}\end{align}
By summing up (\ref{ineq_genp1}) for $k=1,2,\ldots,m$ with an arbitrary
natural number $m\le N$, it holds that
\begin{align}\label{ineq_branch}\begin{aligned}
    &\quad\, \|v_m\|_0^2-\|v_0\|_0^2
    +\tau\sum^m_{k=1}\left(
        \frac{\tau}{4}\|D_\tau v_k\|_0^2
        + \frac{1}{2\tau}\|\vpre_k-v_k\|_0^2
      + \frac{1}{c_a}\|\vpre_k\|_1^2\right)\\
    &\le \tau\sum^m_{k=1}\left(
      c_a\|F_k\|_{H^*}^2
    + 6\tau\|\nabla Q_k\|_0^2
    + 2\langle G_k, \vpre_k \rangle_H \right).
\end{aligned}\end{align}
From the assumption (\ref{asmp_G}) with $\delta:=\frac{1}{4c_a}$;
\[
    \langle G_k, \vpre_k \rangle_H
    \le \frac{\|\vpre_k\|_1^2}{4c_a}
    + A_{\frac{1}{4c_a}}(2\|v_{k-1}\|_0^2 
    + 2\|v_{k-1}-v^*_{k-1}\|_0^2 + \beta_k^2),
\]
we obtain 
\[\begin{aligned}
    & \|v_m\|_0^2 - \|v_0\|_0^2
    +\tau\sum^m_{k=1}\left(
        \frac{\tau}{4}\|D_\tau v_k\|_0^2
        + \frac{1-8\tau A_{\frac{1}{4c_a}}}{2\tau}\|v_k-\vpre_k\|_0^2
        + \frac{1}{2c_a}\|\vpre_k\|_1^2\right)\\
    \le& \tau\sum^{m-1}_{k=0}4A_{\frac{1}{4c_a}}\|v_k\|_0^2
    + \tau\sum^m_{k=1}
    \left(c_a\|F_k\|_{H^*}^2 + 6\tau\|\nabla Q_k\|_0^2
    + 2A_{\frac{1}{4c_a}}\beta_k^2\right),
\end{aligned}\]
where we have used $v_0-\vpre_0=0$.
By the discrete Gronwall inequality, 
if $\tau\le\tau_0:=1/(16A_{\frac{1}{4c_a}})$, 
then it holds that for all $m=0,1,\ldots,N,$
\[\begin{array}{rl}
    &{\DS \|v_m\|_0^2
    +\tau\sum^m_{k=1}\left(
        \frac{\tau}{4}\|D_\tau v_k\|_0^2
        + \frac{1}{4\tau}\|v_k-\vpre_k\|_0^2
        + \frac{1}{2c_a}\|\vpre_k\|_1^2\right)}\\
    \le&{\DS \exp\left(\frac{16}{3}A_{\frac{1}{4c_a}}\right)
    \left\{\|v_0\|_0^2
        + \tau\sum^m_{k=1}
        \left(c_a\|F_k\|_{H^*}^2 + 6\tau\|\nabla Q_k\|_0^2
        + 2A_{\frac{1}{4c_a}}\beta_k^2\right)\right\}},
\end{array}\]
which implies that
\[\begin{array}{rl}
    &{\DS \|\bar{v}_\tau(t)\|_0^2
    +\int^t_0 \left(
        \tau\left\|\frac{\partial \hat{v}_\tau}{\partial t}(s)\right\|_0^2
        + \frac{1}{\tau}\|\bar{v}_\tau(s)-\bar{v}^*_\tau(s)\|_0^2
        + \|\bar{v}^*_\tau(s)\|_1^2\right)ds}\\
    \le&{\DS 
    c_1\left\{\|v_0\|_0^2
    +\int^t_0 (\|\bar{F}_\tau(s)\|_{H^*}^2+\tau\|\bar{Q}_\tau(s)\|_1^2
    +\bar{\beta}_\tau^2(s))ds\right\}}
\end{array}\]
for all $t\in(0,T]$, where 
$c_1:=\exp(16A_{\frac{1}{4c_a}}/3)
\times\max\{c_a,6,2A_{\frac{1}{4c_a}}\}\times\max\{4,2c_a\}$. Hence, 
\begin{align}\label{ineq_genp2}
    \DS \|\bar{v}_\tau\|_{L^\infty(L^2)}^2
    \le M,\qquad
    \tau\left\|\frac{\partial \hat{v}_\tau}{\partial t}\right\|_{L^2(L^2)}^2
    + \frac{1}{\tau}\|\bar{v}_\tau - \bar{v}^*_\tau\|_{L^2(L^2)}^2
    + \|\bar{v}^*_\tau\|_{L^2(H^1)}^2
    \le M,
\end{align}
where 
$
    M := c_1(\|v_0\|_0^2+\|\bar{F}_\tau\|_{L^2(H^*)}^2
        +\tau\|\bar{Q}_\tau\|_{L^2(H^1)}^2
        +\|\bar{\beta}_\tau\|_{L^2(0,T)}^2).
$
If $\langle G_k, \vpre_k \rangle_H$ $\le 0$ for all $k = 1,2,\ldots,N$, 
then we immediately obtain (\ref{ineq_gen}) for all $0<\tau<T$ 
from (\ref{ineq_branch}).

Since it holds that for all $m=1,2,\ldots,N$,
\[
    \|\vpre_m\|_0^2
    \le 2(\|v_m\|_0^2+\|v_m-\vpre_m\|_0^2)
    \le 2\max_{k=1,\ldots,N}\|v_k\|_0^2
    +\tau\sum^N_{k=1}\frac{2}{\tau}\|v_k-\vpre_k\|_0^2,
\]
we obtain for all $0<\tau<\tau_0$,
\[
    \|\bar{v}^*_\tau\|_{L^\infty(L^2)}^2
    \le 2\|\bar{v}_\tau\|_{L^\infty(L^2)}^2
        +\frac{2}{\tau}\|\bar{v}_\tau-\bar{v}^*_\tau\|_{L^2(L^2)}^2
    \le 4M.
\]
\qed
\medskip

By using Lemma \ref{lem_gen}, we prove Theorem \ref{thm_stab}.

\noindent{\it Proof of Theorem \ref{thm_stab}.}
We set $(F_k)_{k=1}^N, (G_k)_{k=1}^N \subset H^*$ defined by 
\[
    \langle F_k, \varphi \rangle_H
    := \langle f_k, \varphi \rangle_H
    - (\nabla p^b_k,\varphi),\qquad
    \langle G_k, \varphi \rangle_H
    := -d(\upre_{k-1}, \upre_k, \varphi)
\]
for all $k=1,2,\ldots,N$ and $\varphi\in H$.
From Problem \ref{def_FS} and the condition (\ref{cond}),
if we set $q_k := P_k - p^b_k$, then 
$(\upre_k, u_k, q_k)_{k=1}^N \subset H \times \Lo^d \times \Hsecond$
satisfies that for all $k = 1,2,\ldots, N$,
\[\left\{\begin{aligned}
    &\frac{1}{\tau}(u_k-u_{k-1},\varphi)
    + a(\upre_k,\varphi) - (q_k,\dvg\varphi)
    = \langle F_k+G_k,\varphi\rangle,\\
    &\tau(\nabla q_k,\nabla\psi) + (\dvg\upre_k,\psi)
    = -(\nabla p^b_k,\nabla \psi),\\
    &u_k = \upre_k - \tau\nabla(q_k+p^b_k) \mbox{ in }\Lo^d,
\end{aligned}\right.\]
with $u_0 \in L^{p_d}(\Omega)^d (\subset \Lo^d)$.
By Lemma \ref{lem_outflow}, it holds that
\[
    \langle G_k, \upre_k \rangle_H
    = -d(\upre_{k-1}, \upre_k, \upre_k) = 0
    \qquad\text{for all }k=1,2,\ldots,N.
\]
Therefore, by Lemma \ref{lem_gen} and Remark \ref{rem_fpb}, 
we conclude the proof.
\qed

\subsection{Convergence}\label{sec_conv}

In this section, we assume Hypothesis \ref{hyp_reg_sol}.
We calculate the error estimates in suitable norms 
between the solutions to (\ref{weak_original}) and (\ref{weak_FS}).
By Hypothesis \ref{hyp_reg_sol} and the first equation of (\ref{weak_original}),
it holds that $\frac{\partial u}{\partial t}\in C([0,T];H^*)$
and, for all $\varphi\in H$ and $k=1,2,\ldots,N$,
\begin{align*}\begin{aligned}
    \frac{1}{\tau}(u(t_k)-u(t_{k-1}),\varphi)
    + a(u(t_k),\varphi) + d(\upre_{k-1},\upre_k,\varphi)
    + (\nabla P(t_k),\varphi)
    = \langle f(t_k) - R_k - R^{\rm n.l.}_k,\varphi\rangle_H,
\end{aligned}\end{align*}
where $R_k, R^{\rm n.l.}_k\in H^*$ defined by
\[\begin{aligned}
    \langle R_k, \varphi \rangle_H
    &:= \left\langle \frac{\partial u}{\partial t}(t_k)
    - \frac{u(t_k)-u(t_{k-1})}{\tau}, \varphi\right\rangle_H,\\
    \langle R^{\rm n.l.}_k, \varphi \rangle_H
    &:= d(u(t_k), u(t_k), \varphi) - d(\upre_{k-1}, \upre_k, \varphi)
\end{aligned}\]
for all $\varphi \in H$.
If we put $e_0=0$, $e_k:=u_k-u(t_k)\in\Lo^d$, 
$e^*_k:=\upre_k-u(t_k)\in H$, and $q_k:=P_k-P(t_k)\in\Hsecond$ 
for $k=1,2,\ldots,N$, by (\ref{weak_FS_equiv2}), then it holds that 
for all $k=1,2,\ldots,N$, $\varphi\in H$, and $\psi\in\Hsecond$,
\begin{align}\label{eq_error}\left\{\begin{aligned}
    &\frac{1}{\tau}(e_k-e_{k-1},\varphi)
    + a(e^*_k,\varphi) - (q_k,\dvg\varphi)
    = \langle R_k + R^{\rm n.l.}_k, \varphi \rangle_H\\
    &\tau(\nabla q_k,\nabla\psi) + (\dvg e^*_k,\psi)
    = -\tau(\nabla P(t_k)\nabla\psi),\\
    &e_k = e^*_k - \tau\nabla (q_k+P(t_k)) \mbox{ in }\Lo^d,
\end{aligned}\right.\end{align}
where we have used $(\nabla q_k,\varphi) = -(q_k,\dvg\varphi)$.

In order to prove Theorems \ref{thm_conv1} and \ref{thm_conv2}, 
we prepare Lemmas \ref{lem_R} and \ref{lem_diffdt}.
See the Appendix for the proofs.

\begin{Lem}\label{lem_R}

  (i) Under Hypothesis \ref{hyp_reg_sol}, we have 
  \[
      \|\bar{R}_\tau\|_{L^2(H^*)}^2
      \le \frac{\tau^2}{3}
      \left\|\frac{\partial^2 u}{\partial t^2}\right\|_{L^2(H^*)}^2.
  \]

  (ii) Furthermore, if Hypothesis \ref{hyp_reg_sol2} holds,
  then we have 
  \[
      \sum^N_{k=2}\tau\|D_\tau R_k\|_{H^*}^2
      \le \frac{2}{3}\tau^2
      \left\|\frac{\partial^3 u}{\partial t^3}\right\|_{L^2(H^*)}^2.
  \]
\end{Lem}

\medskip

\begin{Lem}\label{lem_diffdt}
    Let $(E,(\cdot,\cdot)_E)$ be a Hilbert space and let $x\in C([0,T];E)$
    satisfy that $\frac{\partial x}{\partial t}\in$\\ $L^2(0,T;E)$.

    \noindent(i) It holds that for all $k=1,2,\ldots,N,$
    \[
        \|D_\tau x(t_k)\|_E
        \le \frac{1}{\sqrt{\tau}}\left\|\frac{\partial x}{\partial t}\right\|_{L^2(t_{k-1},t_k;E)}.
    \]
    (ii) It holds that 
    \[
        \|x-x_\tau\|_{L^\infty(E)}
        \le\sqrt{\tau}\left\|\frac{\partial x}{\partial t}\right\|_{L^2(E)},\qquad
        \|x-x_\tau\|_{L^2(E)}
        \le\frac{\tau}{\sqrt{2}}\left\|\frac{\partial x}{\partial t}\right\|_{L^2(E)}.
    \]
\end{Lem}

By using Lemmas \ref{lem_gen}, \ref{lem_R} and \ref{lem_diffdt}, 
we prove Theorem \ref{thm_conv1}.

\noindent{\it Proof of Theorem \ref{thm_conv1}.}
For all $\delta>0$ and $k=1,2,\ldots,N$, 
by Lemmas \ref{lem_outflow}, \ref{lem_nonlin} and \ref{lem_diffdt}, 
we have 
\begin{align}\label{nl_sepa}\begin{aligned}
  \langle R^{\rm n.l.}_k, e^*_k \rangle_H
  &= -d(\upre_{k-1},e^*_k,e^*_k)
  - d(e^*_{k-1},u(t_k),e^*_k)
  + \tau d(D_\tau u(t_k),u(t_k),e^*_k)\\
  &\le c_d\|e^*_{k-1}\|_0\|u(t_k)\|_2\|e^*_k\|_1
  + c_d\tau\|D_\tau u(t_k)\|_0\|u(t_k)\|_2\|e^*_k\|_1\\
  &\le \frac{\delta}{2}\|e^*_k\|_1^2
  + \frac{c_d^2\|u(t_k)\|_2^2}{2\delta}\|e^*_{k-1}\|_0^2
  + \frac{\delta}{2}\|e^*_k\|_1^2
  + \frac{c_d^2\|u(t_k)\|_2^2}{2\delta}\tau^2\|D_\tau u(t_k)\|_0^2\\
  &\le \delta\|e^*_k\|_1^2
  + \frac{c_d^2 c_{\rm max}^2}{2\delta}\|e^*_{k-1}\|_0^2
  + \frac{c_d^2c_{\rm max}^2}{2\delta}\tau\left\|\frac{\partial u}{\partial t}\right\|_{L^2(t_{k-1},t_k;\Lo^d)}^2
\end{aligned}\end{align} 
where $c_{\rm max}:=\|u\|_{C([0,T],H^2(\Omega)^d)}$.
By (\ref{eq_error}) and Lemmas \ref{lem_gen}, \ref{lem_R} and \ref{lem_diffdt}, 
there exist two constants $\tau_0,c_1>0$ such that
for all $0<\tau<\tau_0$,
\[\begin{aligned}
    & \|\bar{e}_\tau\|_{L^\infty(L^2)}^2
    + \|\bar{e}^*_\tau\|_{L^\infty(L^2)}^2
    + \|\bar{e}^*_\tau\|_{L^2(H^1)}^2
    + \frac{1}{\tau}\|\bar{e}_\tau-\bar{e}^*_\tau\|_{L^2(L^2)}^2\\
    \le\,& c_1\Biggl(\|\bar{R}_\tau\|_{L^2(H^*)}^2
    + \tau\|P_\tau\|_{L^2(H^1)}^2
    + \tau^2\left\|\frac{\partial u}{\partial t}\right\|_{L^2(L^2)}^2\Biggr)\\
    \le\,& c_1\left(\frac{\tau^2}{3}\left\|\frac{\partial^2 u}{\partial t^2}\right\|_{L^2(H^*)}^2
    + 2\tau\|P\|_{L^2(H^1)}^2
    + \tau^3\left\|\frac{\partial P}{\partial t}\right\|_{L^2(H^1)}^2
    + \tau^2\left\|\frac{\partial u}{\partial t}\right\|_{L^2(L^2)}^2\right),
\end{aligned}\]
which implies that
\[\begin{aligned}
    \|\bar{u}_\tau - u_\tau\|_{L^\infty(L^2)}
    +\|\bar{u}^*_\tau - u_\tau\|_{L^\infty(L^2)}
    +\|\bar{u}^*_\tau - u_\tau\|_{L^2(H^1)}
    &\le c_2\sqrt{\tau},\\
    \|\bar{u}_\tau-\bar{u}^*_\tau\|_{L^2(L^2)}
    &\le c_2\tau
\end{aligned}\]
for a constant $c_2>0$,
where we have used $\bar{e}_\tau = \bar{u}_\tau- u_\tau$
and $\bar{e}^*_\tau = \bar{u}^*_\tau - u_\tau$.
By the triangle inequality and Lemma \ref{lem_diffdt}, 
it holds that 
$\|u-\bar{u}_\tau\|_{L^\infty(L^2)}
+ \|u-\bar{u}^*_\tau\|_{L^\infty(L^2)}
\le c_2\sqrt{\tau} 
+ 2\sqrt{\tau}\|\frac{\partial u}{\partial t}\|_{L^2(L^2)}$.
To complete the first inequality of Theorem \ref{thm_conv1},
it is sufficient to prove that $\|u-u_\tau\|_{L^2(H^1)}\le c_3\sqrt{\tau}$
for a constant $c_3>0$.
Since $u(t)\in H\cap H^2(\Omega)^d$ and $\dvg u(t)=0\in\Hb$ 
for all $t\in[0,T]$, by Proposition \ref{prop_weak},
Lemmas \ref{lem_a} and \ref{lem_diffdt}, we find that 
\[\begin{aligned}
    &\quad \|u-u_\tau\|_{L^2(H^1)}^2
    = \sum^N_{k=1}\int^{t_k}_{t_{k-1}}\|u(t)-u(t_k)\|_1^2 dt\\
    &\le c_a\sum^N_{k=1}\int^{t_k}_{t_{k-1}}
    a(u(t)-u(t_k), u(t)-u(t_k)) dt 
    = c_a\sum^N_{k=1}\int^{t_k}_{t_{k-1}}
    \left(-\Delta(u(t)-u(t_k)), u(t)-u(t_k)\right) dt\\
    &\le 2\sqrt{d}c_a c_{\rm max}\int^T_0\|u(t)-u(t_k)\|_0 dt 
    \le 2\sqrt{dT}c_a c_{\rm max}\|u-u_\tau\|_{L^2(L^2)} 
    \le \sqrt{2dT}c_a c_{\rm max}\tau
    \left\|\frac{\partial u}{\partial t}\right\|_{L^2(L^2)}.
\end{aligned}\]
\qed
\medskip

We improve the error estimates for the velocity and 
pressure in the $L^2(L^2)$-norm.
In order to prove Theorem \ref{thm_conv2},
we prepare Proposition \ref{prop_stokes_reg} and Lemma \ref{lem_AN}.

\begin{Prop}\label{prop_stokes_reg}
    Under Hypothesis \ref{hyp_stokes}, for all $e\in\Lo$, 
    the pair of functions $(w,r)=T(e)$
    belongs to $H^2(\Omega)\times\Hsecond$.
\end{Prop}
\noindent{\it Proof.}
By Hypothesis \ref{hyp_stokes}, $(w,r)\in H^2(\Omega)^d\times\Hb$.
Since it holds that for all $\varphi\in H$,
\[
    0=a(w,\varphi)-(r,\dvg\varphi)-(e,\varphi)
    =\into(\nabla\times(\nabla\times w)+\nabla r-e)\cdot\varphi dx
    -\int_{\Gamma_2}r\varphi\cdot n ds,
\]
we obtain $r\in\Hsecond$.
\qed

\begin{Lem}\label{lem_AN}
    Under the assumption of Lemma \ref{lem_gen} and Hypothesis \ref{hyp_stokes},
    if we assume the following conditions: 
    if $(w_k,r_k):=T(v_k)$ for all $k=0,1,\ldots,N$, then 
    for all $\delta>0$ there exist a constant $A_\delta>0$
    independent of $k$ and $\tau$, and a sequence $(\gamma_k)_{k=1}^N\in\R$ 
    such that for all $k=1,2,\ldots,N$,
    \begin{align}\label{asmp_G2}
        \langle G_k, w_k \rangle_H
        \le \delta(\|v^*_{k-1}\|_0^2+\|v^*_k\|_0^2)
        +A_\delta(\|w_k\|_1^2+\gamma_k^2),
    \end{align}
    \noindent then there exist two constants 
    $\tau_0,c>0$ independent of $\tau$ such that for all $0<\tau<\tau_0$,
    \[
        \|\bar{v}_\tau\|_{L^2(L^2)}^2
        \le c(\|v_0\|_{V^*}^2+\tau\|\vpre_0\|_0^2
        +\|\bar{v}^*_\tau-\bar{v}_\tau\|_{L^2(L^2)}^2
        +\|\bar{F}_\tau\|_{L^2(H^*)}^2+\|\bar{\gamma}_\tau\|_{L^2(0,T)}^2).
    \]
\end{Lem}
\noindent{\it Proof.}
Let $(w_k,r_k):=T(v_k)$ for all $k=0,1,\ldots,N$.
It follows from Proposition \ref{prop_stokes_reg} that $r_k \in \Hsecond$.
The first equation of (\ref{eq_gen}) implies that 
for all $k=1,2,\ldots,N$,
\begin{align}\label{eq_gen13}
    \frac{1}{\tau}(v_k-v_{k-1},w_k) + a(\vpre_k,w_k) 
    = \langle F_k, w_k \rangle_H + \langle G_k, w_k \rangle_H,
\end{align}
where we have used $\dvg w_k=0$ in $\Lo$.
By Lemma \ref{lem_stokes}, we obtain
\[\begin{aligned}
    (v_k-v_{k-1},w_k)
    &=a(w_k,w_k)-(r_k,\dvg w_k)-a(w_{k-1},w_k)+(r_{k-1},\dvg w_k)\\
    &=\frac{1}{2}\left(\|w_k\|_a^2 - \|w_{k-1}\|_a^2
        + \|w_k-w_{k-1}\|_a^2\right)
    \ge \frac{c_1}{2}\left(\|w_k\|_1^2-\|w_{k-1}\|_1^2\right)
\end{aligned}\]
where $c_1 := \min\{c_a, c_a^{-1}\}$.
For the second term of the left hand side of (\ref{eq_gen13}), 
by the definition of the operator $T$, we have
\[\begin{aligned}
    a(\vpre_k,w_k)
    = \|v_k\|_0^2+(v_k,\vpre_k-v_k)-(\nabla r_k,\vpre_k-v_k),
\end{aligned}\]
where we have used 
the third equation of (\ref{eq_gen}) and
$(\nabla r_k,v_k)=(\nabla r_k,\vpre_k)-\tau(\nabla r_k,\nabla (q_k+Q_k))=0$.
By Hypothesis \ref{hyp_stokes}, it holds that 
\[\begin{aligned}
    |(v_k,\vpre_k-v_k)-(\nabla r_k,\vpre_k-v_k)|
    &\le (\|v_k\|_0+\|\nabla r_k\|_0)\|\vpre_k-v_k\|_0
    \le c_2\|v_k\|_0\|\vpre_k-v_k\|_0\\
    &\le \frac{1}{4}\|v_k\|_0^2+c_2^2\|\vpre_k-v_k\|_0^2
\end{aligned}\]
for a constant $c_2>0$. Hence, we have
\[
    a(\vpre_k,w_k)
    \ge \frac{3}{4}\|v_k\|_0^2-c_2^2\|\vpre_k-v_k\|_0^2.
\]
For the first term of the right hand side of (\ref{eq_gen13}), 
by Lemma \ref{lem_stokes}, we have
\[
    \langle F_k,w_k\rangle_H
    \le \|F_k\|_{H^*}\|w_k\|_1
    \le c_3\|F_k\|_{H^*}\|v_k\|_0
    \le \frac{1}{4}\|v_k\|_0^2+c_3^2\|F_k\|_{H^*}^2
\]
for a constant $c_3>0$.
Hence, we have that for all $k=1,2,\ldots,N$,
\[
    \|w_k\|_1^2 - \|w_{k-1}\|_1^2 + \frac{\tau}{c_1}\|v_k\|_0^2
    \le \frac{2\tau}{c_1}(c_2^2 \|\vpre_k-v_k\|_0^2 
    + c_3^2 \|F_k\|_{H^*}^2
    + \langle G_k,w_k\rangle_H).
\]
By summing up for $k=1,2,\ldots,m$ with an arbitrary natural number $m\le N$,
it holds that
\begin{align*}\begin{aligned}
    \|w_m\|_1^2 \!-\! \|w_0\|_1^2
    \!+\! \frac{\tau}{c_1}\sum^m_{k=1}\|v_k\|_0^2
    \!\le\! \frac{2\tau}{c_1}\sum^m_{k=1}(c_2^2 \|\vpre_k-v_k\|_0^2 
    \!+\! c_3^2 \|F_k\|_{H^*}^2
    \!+\! \langle G_k,w_k\rangle_H).
\end{aligned}\end{align*}
From the assumption (\ref{asmp_G2}) with $\delta:=\frac{1}{16}$,
we obtain for all $m=1,2,\ldots,N$,
\[\begin{aligned}
    \sum^m_{k=1}\langle G_k,w_k\rangle_H
    &\le \sum^m_{k=1}\left\{\frac{1}{16}(\|v^*_{k-1}\|_0^2+\|v^*_k\|_0^2)
    + A_{\frac{1}{16}}(\|w_k\|_1^2 + \gamma_k^2)\right\}\\
    &\le \frac{1}{16}\|\vpre_0\|_0^2
    + \frac{1}{4}\sum^m_{k=1}(\|v_k\|_0^2 + \|v^*_k-v_k\|_0^2)
    + A_{\frac{1}{16}}\sum^m_{k=1}(\|w_k\|_1^2 + \gamma_k^2)
\end{aligned}\]
and hence, 
\[\begin{aligned}
    &\quad \|w_m\|_1^2-\|w_0\|_1^2
    + \frac{\tau}{2c_1}\sum^m_{k=1}\|v_k\|_0^2\\
    &\le \tau\sum^m_{k=1}\frac{2A_{\frac{1}{16}}}{c_1} \|w_k\|_1^2
    + \frac{\tau}{8c_1}\|\vpre_0\|_0^2
    + \tau\sum^m_{k=1}c_4(\|\vpre_k-v_k\|_0^2
    + \|F_k\|_{H^*}^2 + \gamma_k^2),
\end{aligned}\]
where $c_4 := 2c_1^{-1}\max\{ c_2^2+1/4, c_3^2, A_{\frac{1}{16}} \}$.
By the discrete Gronwall inequality, 
if $\tau\le\tau_0 := c_1/A_{\frac{1}{16}}$,
then we have 
\[\begin{aligned}
    &\quad \|w_N\|_1^2 + \frac{\tau}{2c_1}\sum^N_{k=1}\|v_k\|^2\\
    &\le \exp\left(\frac{4A_{\frac{1}{16}}}{c_1}\right)
    \biggl\{\|w_0\|_1^2 + \frac{\tau}{8c_1}\|\vpre_0\|_0^2
    + \tau\sum^N_{k=1}c_4\left(\|\vpre_k-v_k\|_0^2
    \!+ \|F_k\|_{H^*}^2 \!+ \gamma_k^2\right)\biggl\}.
\end{aligned}\]
Therefore, by Lemma \ref{lem_stokes}, we obtain
\[\begin{aligned}
    \|\bar{v}_\tau\|_{L^2(L^2)}^2
    \le c_5\Bigl(\|v_0\|_{V^*}^2+\tau\|\vpre_0\|_0^2
    +\|\bar{v}^*_\tau-\bar{v}_\tau\|_{L^2(L^2)}^2
    +\|F_\tau\|_{L^2(H^*)}^2+\|\bar{\gamma}_\tau\|_{L^2}^2\Bigr)
\end{aligned}\]
for a constant $c_5>0$.
\qed
\medskip

We prove the first inequality of Theorem \ref{thm_conv2}

\noindent{\it Proof of the first inequality of Theorem \ref{thm_conv2}.}
We apply Lemmas \ref{lem_AN} for (\ref{eq_error}).
Let $(w_k,r_k) := T(e_k)$ for all $k = 0,1,\ldots,N$.
It holds that for all $k=1,2,\ldots,N$,
\[
    \langle R^{\rm n.l.}_k, w_k \rangle_H
    = -d(e^*_{k-1},\upre_k,w_k) - d(u(t_{k-1}),e^*_k,w_k)
    + \tau d(D_\tau u(t_k),u(t_k),w_k).
\]
Hypothesis \ref{hyp_stokes} and Theorem \ref{thm_conv1} implies that
there exists a constant $c_1>0$ such that 
$\|w_k\|_2\le c_1\sqrt{\tau}$ for all $k=1,2,\ldots,N$.
It holds that for all $\delta>0$,
\[\begin{aligned}
    -d(e^*_{k-1},\upre_k,w_k)
    =\,& -d(e^*_{k-1},u(t_k),w_k) - d(e^*_{k-1},e^*_k,w_k)\\
    \le\,& c_d\|e^*_{k-1}\|_0\|u(t_k)\|_2\|w_k\|_1
    + c_d\|e^*_{k-1}\|_0\|e^*_k\|_1\|w_k\|_2\\
    \le\,& c_d c_{\rm max}\|e^*_{k-1}\|_0\|w_k\|_1
    + c_d c_1\sqrt{\tau}\|e^*_{k-1}\|_0\|e^*_k\|_1\\
    \le\,& \frac{\delta}{2}\|e^*_{k-1}\|_0^2
    + \frac{c_d^2 c_{\rm max}^2}{2\delta}\|w_k\|_1^2
    + \frac{\delta}{2}\|e^*_{k-1}\|_0^2
    + \frac{c_d^2 c_1^2}{2\delta}\tau\|e^*_k\|_1^2\\
    \le\,& \delta\|e^*_{k-1}\|_0^2
    + \frac{c_d^2 c_{\rm max}^2}{2\delta}\|w_k\|_1^2
    + \frac{c_d^2 c_1^2}{2\delta}\tau\|e^*_k\|_1^2,\\
    -d(u(t_{k-1}),e^*_k,w_k)
    =\,& d(u(t_{k-1}),w_k,e^*_k)
    \le c_d\|u(t_{k-1})\|_2\|w_k\|_1\|e^*_k\|_0\\
    \le\,& \delta\|e^*_k\|_0^2
    + \frac{c_d^2 c_{\rm max}^2}{4\delta}\|w_k\|_1^2\\
    \tau d(D_\tau u(t_k),u(t_k),w_k)
    \le\,& c_d\tau\|D_\tau u(t_k)\|_0\|u(t_k)\|_2\|w_k\|_1
    \le c_d c_{\rm max}\tau\|D_\tau u(t_k)\|_0\|w_k\|_1\\
    \le\,& \frac{c_d^2 c_{\rm max}^2}{4\delta}\|w_k\|_1^2
    + \delta\tau^2\|D_\tau u(t_k)\|_0^2,
\end{aligned}\]
where $c_{\rm max}:=\|u\|_{C([0,T];H^2(\Omega)^d)}$.
Hence, by Lemma \ref{lem_diffdt}, it holds that for all $k=1,2,\ldots,N$,
\[\begin{aligned}
    \langle R^{\rm n.l.}_k,w_k\rangle_H
    &\le \delta\|e^*_{k-1}\|_0^2+\delta\|e^*_k\|_0^2
    + c_\delta(\|w_k\|_1^2 + \tau\|e^*_k\|_1^2 + \tau^2\|D_\tau u(t_k)\|_0^2)\\
    &\le \delta\|e^*_{k-1}\|_0^2 + \delta\|e^*_k\|_0^2
    + c_\delta\left(\|w_k\|_1^2 + \tau\|e^*_k\|_1^2
    + \tau\left\|\frac{\partial u}{\partial t}\right\|_{L^2(t_{k-1},t_k;\Lo^d)}^2\right),
\end{aligned}\]
where 
$c_\delta:=\max\{ c_d^2 c_{\rm max}^2 / \delta,
    c_d^2 c_1^2 / (2\delta), \delta\}$.
By Lemma \ref{lem_AN}, there exist two constants $c_2,\tau_0>0$ 
such that for all $0<\tau<\tau_0$,
\[\begin{aligned}
    \|\bar{e}_\tau\|_{L^2(L^2)}^2
    \le c_2\left(\|\bar{e}^*_\tau-\bar{e}_\tau\|_{L^2(L^2)}^2
    + \|\bar{R}_\tau\|_{L^2(H^*)}^2
    + \tau\|\bar{e}^*_\tau\|_{L^2(H^1)}^2
    + \tau^2\left\|\frac{\partial u}{\partial t}\right\|_{L^2(L^2)}^2\right).
\end{aligned}\]
By Lemma \ref{lem_R} and Theorem \ref{thm_conv1},
there exists a constant $c_3>0$ such that for all $0<\tau<\tau_0$,
\[
    \|u_\tau-\bar{u}_\tau\|_{L^2(L^2)} \le c_3\tau.
\]
By Lemma \ref{lem_diffdt} and Theorem \ref{thm_conv1}, 
we obtain the first inequality of Theorem \ref{thm_conv2};
\[\begin{aligned}
    \|u-\bar{u}_\tau\|_{L^2(L^2)}+\|u-\bar{u}^*_\tau\|_{L^2(L^2)}
    \le 2\|u-u_\tau\|_{L^2(L^2)}+2\|u_\tau-\bar{u}_\tau\|_{L^2(L^2)}
    +\|\bar{u}_\tau-\bar{u}^*_\tau\|_{L^2(L^2)}
    \le c_4\tau
\end{aligned}\]
for a constant $c_4>0$.
\qed
\medskip

To prove the second inequality of Theorem \ref{thm_conv2}, 
we prepare the following two lemmas:

\begin{Lem}\label{lem_firstDt}
    Under Hypothesis \ref{hyp_reg_sol}, there exists a constant $c>0$
    independent of $\tau$ such that 
    \[
        \|D_\tau e_1\|_{V^*}\le c\sqrt{\tau},\qquad
        \|D_\tau e_1\|_0+\|D_\tau e^*_1\|_0\le c,\qquad
        \|D_\tau e^*_1\|_1\le \frac{c}{\sqrt{\tau}}.
    \]
\end{Lem}

\noindent{\it Proof.}
By (\ref{eq_error}) and (\ref{ineq_genp1}) with $k:=1$
in the proof of Lemma \ref{lem_gen}, we obtain
\[\begin{aligned}
    \|e_1\|_0^2
    \!+\! \frac{1}{2}\|e_1-e^*_1\|_0^2 \!+\! \frac{\tau}{c_a}\|e^*_1\|_1^2
    \!\le\! c_a\tau\|R_1\|_{H^*}^2
    \!+\! 6\tau^2\|\nabla P(t_1)\|_0^2
    \!+\! 2\tau\langle R^{\rm n.l.}_1,e^*_1\rangle_H.
\end{aligned}\]
Putting $k:=1$ and $\delta:=\frac{1}{4c_a}$ in (\ref{nl_sepa}), 
it holds that 
\[
    \langle R^{\rm n.l.}_1,e^*_1\rangle_H
    \le \frac{1}{4c_a}\|e^*_1\|_1^2
    + 2c_a c_d^2 c_{\rm max}^2\tau
    \left\|\frac{\partial u}{\partial t}\right\|_{L^2(0,t_1;\Lo^d)}^2
\]
where $c_{\rm max}:=\|u\|_{C([0,T];H^2(\Omega)^d)}$.
Hence, by Lemma \ref{lem_R}, we have
\[\begin{aligned}
    &\quad \|e_1\|_0^2 + \frac{1}{2}\|e_1-e^*_1\|_0^2
    + \frac{1}{2c_a}\tau\|e^*_1\|_1^2\\
    &\le c_a\tau\|R_1\|_{H^*}^2
    + 6\tau^2\|\nabla P(t_1)\|_0^2
    + 4c_a c_d^2 c_{\rm max}^2\tau^2\left\|\frac{\partial u}{\partial t}\right\|_{L^2(0,t_1;\Lo^d)}^2\\
    &\le \frac{c_a \tau^2}{3}\left\|\frac{\partial u}{\partial t}\right\|_{L^2(H^*)}^2
    + 6\tau^2\|P\|_{C([0,T];H^1)}^2
    + 4c_a c_d^2 c_{\rm max}^2\tau^2\left\|\frac{\partial u}{\partial t}\right\|_{L^2(L^2)}^2
    \le c_2\tau^2
\end{aligned}\]
where $c_2:=c_a(\frac{1}{3} + 4c_d^2 c_{\rm max}^2)\left\|\frac{\partial u}{\partial t}\right\|_{L^2(L^2)}^2
+ 6\|P\|_{C([0,T];H^1)}^2$, which implies that 
$\|D_\tau e_1\|_0 =$\\ $\tau^{-1}\|e_1\|_0 \le \sqrt{c_2}$,
$\|D_\tau e^*_1\|_1 \le \sqrt{2c_a c_2}\tau^{-1/2}$ and 
\[
    \|D_\tau e^*_1\|_0
    = \frac{1}{\tau}\|e^*_1\|_0
    \le \frac{1}{\tau}(\|e_1\|_0 + \|e_1-e^*_1\|_0)
    \le (1+\sqrt{2})\sqrt{c_2}.
\]
On the other hand, by (\ref{eq_error}) and 
Lemmas \ref{lem_a}, \ref{lem_R},
\begin{align*}
    \|D_\tau e_1\|_{V^*}
    &= \sup_{0\ne\varphi\in V}\frac{|(e_1-e_0,\varphi)|}{\tau\|\varphi\|_1}\\
    &=\!\! \sup_{0\ne\varphi\in V}\!\frac{\left|
    - a(e^*_1,\varphi) \!+\! (q_1,\dvg\varphi)
    \!+\! \langle R_1, \!\varphi\rangle_H
    \!-\! d(u_0,\!e^*_1,\!\varphi) \!+\! \tau d(D_\tau u(t_1),\!u(t_1),\!\varphi)\right|}{\|\varphi\|_1}\\
    &\le c_a\|e^*_1\|_1 + \|R_1\|_{H^*} 
    + c_d(\|u_0\|_1\|e^*_1\|_1 + \tau\|D_\tau u(t_1)\|_0\|u(t_1)\|_2)\\
    &\le c_a\|e^*_1\|_1  
    \!+\! \sqrt{\frac{\tau}{3}}\left\|\frac{\partial u}{\partial t}\right\|_{L^2(0,t_1;H^*)}
    \!+\! c_d c_{\rm max}\left(\|e^*_1\|_1
    + \sqrt{\tau} \left\|\frac{\partial u}{\partial t}\right\|_{L^2(0,t_1,\Lo^d)}\right)\\
    &\le \sqrt{\tau}\left\{(c_a+c_d c_{\rm max})\sqrt{2c_a c_2}
    + \left(\frac{1}{\sqrt{3}} + c_d c_{\rm max}\right)
    \left\|\frac{\partial u}{\partial t}\right\|_{L^2(L^2)}\right\},
\end{align*}
where $c_{\rm max}:=\|u\|_{C([0,T];H^2(\Omega)^d)}$.
\qed
\medskip

\begin{Lem}\label{lem_dedt}
  Under Hypothesis \ref{hyp_reg_sol}, \ref{hyp_stokes}, 
  and \ref{hyp_reg_sol2},
  there exist two constants $c,\tau_0>0$ independent of $\tau$ 
  such that for all $0<\tau<\tau_0$, 
  \[
      \left\|\frac{\partial \hat{e}_\tau}{\partial t}\right\|_{L^2(L^2)}
      \le c\sqrt{\tau}.
  \]
\end{Lem}

\noindent{\it Proof.} 
By (\ref{eq_error}), it holds that 
$(D_\tau e^*_k, D_\tau q_k, D_\tau e_k)_{k=2}^N 
\subset H \times \Hsecond \times \Lo^d$
and for all $k=2,3,\ldots, N$,
$\varphi\in H$, and $\psi\in\Hsecond$,
\begin{align}\label{eq_dt}\left\{\begin{aligned}
  &\!\!\left(\! \frac{D_\tau e_k-D_\tau e_{k-1}}{\tau},\varphi \!\right)
  \!+ a(D_\tau e^*_k,\varphi)
  - (D_\tau q_k,\dvg\varphi)
  \!=\! \langle D_\tau R_k+D_\tau R^{\rm n.l.}_k,\varphi\rangle_H,\\
  &\tau(\nabla D_\tau q_k,\nabla\psi)
  + (\dvg D_\tau e^*_k,\psi)
  = -(\nabla D_\tau P(t_k),\nabla\psi),\\
  &D_\tau e_k=D_\tau e^*_k-\tau\nabla D_\tau (q_k+P(t_k))
  \mbox{ in }\Lo^d
\end{aligned}\right.\end{align}
with $D_\tau e_1 = \tau^{-1}(e_1-e_0) = \tau^{-1}e_1$.
It holds for all $k=2,3,\ldots,N$ and $\varphi\in H$,
\begin{align}\label{eq_dtnl}\begin{aligned}
    \tau\langle D_\tau R^{\rm n.l.}_k,\varphi\rangle_H
    &= -\tau d(\upre_{k-2},D_\tau e^*_k,\varphi)
    - \tau d(D_\tau \upre_{k-1},e^*_k,\varphi)
    + d(e^*_{k-2},u(t_{k-1}),\varphi)\\
    &\quad -d(e^*_{k-1},u(t_k),\varphi)
    - \tau d(D_\tau u(t_{k-1}),u(t_{k-1}),\varphi)
    + \tau d(D_\tau u(t_k),u(t_k),\varphi).
\end{aligned}\end{align}
Here, by Lemma \ref{lem_nonlin}, the right hand side
except for the first and second terms are evaluated from above 
for all $k=2,3,\ldots,N$, $\varphi\in H$ and $\delta>0$,
\begin{align}\label{ineq_diffnl}\begin{aligned}
    &
    d(e^*_{k-2},u(t_{k-1}),\varphi)
    - d(e^*_{k-1},u(t_k),\varphi)
    - \tau d(D_\tau u(t_{k-1}),u(t_{k-1}),\varphi)
    + \tau d(D_\tau u(t_k),u(t_k),\varphi)\\
    \le\,& c_dc_{\rm max}\left(\|e^*_{k-2}\|_0
    +\|e^*_{k-1}\|_0
    +\tau\|D_\tau u(t_{k-1})\|_0
    +\tau\|D_\tau u(t_k)\|_0\right)\|\varphi\|_1\\
    \le\,&\frac{\delta}{2}\|\varphi\|_1^2
    +\frac{2c_d^2 c_{\rm max}^2}{\delta}\sum_{i=0}^1
    \left(\|e^*_{k-i-1}\|_0^2+\tau^2\|D_\tau u(t_{k-i})\|_0^2\right),
\end{aligned}\end{align}
where $c_{\rm max}:=\|u\|_{C([0,T];H^2(\Omega)^d)}$.
By Lemma \ref{lem_outflow}, it holds that 
\[
    -\tau d(\upre_{k-2},D_\tau e^*_k,D_\tau e^*_k)= 0.
\]
By Theorem \ref{thm_conv1}, there exist two constants $\tau_1, c_1>0$
such that $\|\bar{e}^*_\tau\|_{L^2(H^1)}\le c_1$ for all $0<\tau<\tau_1$, 
and hence for all $k=1,2,\ldots, N$, $\|e^*_k\|_1\le c_1$ and
\[\begin{aligned}
    -\tau d(D_\tau u^*_{k-1},e^*_k,D_\tau e^*_k)
    =\,& -\tau d(D_\tau u(t_{k-1}),e^*_k,D_\tau e^*_k)
    - d(D_\tau e^*_{k-1},e^*_k,e^*_k)
    + d(D_\tau e^*_{k-1},e^*_k,e^*_{k-1})\\
    \le\,& c_d\tau \|D_\tau u(t_{k-1})\|_1\|e^*_k\|_1\|D_\tau e^*_k\|_1
    +c_d \|D_\tau e^*_{k-1}\|_1\|e^*_k\|_1\|e^*_{k-1}\|_1\\
    \le\,& c_d c_1\tau \|D_\tau u(t_{k-1})\|_1\|D_\tau e^*_k\|_1
    +c_d c_1\|D_\tau e^*_{k-1}\|_1\|e^*_{k-1}\|_1\\
    \le\,& \frac{\delta}{2}\|D_\tau e^*_k\|_1^2
    +\frac{c_d^2c_1^2}{2\delta}\tau^2\|D_\tau u(t_{k-1})\|_1^2
    +\delta\|D_\tau e^*_{k-1}\|_1^2
    +\frac{c_d^2c_1^2}{4\delta}\|e^*_{k-1}\|_1^2.
\end{aligned}\]
Hence, by (\ref{eq_dtnl}) with $\varphi:=D_\tau e^*_k$
and Lemma \ref{lem_diffdt}, 
for all $0<\tau<\tau_1$, $k=2,3,\ldots,N$ and $\delta>0$,
\[\begin{aligned}
    &\tau\langle D_\tau R^{\rm n.l.}_k,D_\tau e^*_k\rangle_H
    \le \delta(\|D_\tau e^*_k\|_1^2 + \|D_\tau e^*_{k-1}\|_1^2)
    +c_\delta\sum_{i=0}^1
    \left(\|e^*_{k-i-1}\|_1^2+\tau^2\|D_\tau u(t_{k-i})\|_1^2\right)\\
    \le\,& \delta(\|D_\tau e^*_k\|_1^2 + \|D_\tau e^*_{k-1}\|_1^2)
    +c_\delta\sum_{i=0}^1
    \left(\|e^*_{k-i-1}\|_1^2
    +\tau\left\|\frac{\partial u}{\partial t}\right\|_{L^2(t_{k-i-1},t_{k-i};\Hb^d)}^2\right),
\end{aligned}\]
where 
$
    c_\delta := \delta^{-1}(2c_d^2c_{\rm max}^2 + 2^{-1}c_d^2c_1^2).
$
Putting $\delta:=1/(4c_a)$, by (\ref{eq_dt}) and (\ref{ineq_genp1})
in the proof of Lemma \ref{lem_gen}, we have for all 
$0<\tau<\tau_1$ and $k=2,3,\ldots,N$,
\[\begin{aligned}
    &\quad\, \|D_\tau e_k\|_0^2 - \|D_\tau e_{k-1}\|_0^2
    + \frac{1}{2}\|D_\tau e^*_k-D_\tau e_k\|_0^2
    + \frac{\tau}{c_a}\|D_\tau e^*_k\|_1^2\\
    &\le c_a \tau\|D_\tau R_k\|_{H^*}^2
    + 6\tau^2\|\nabla D_\tau P(t_k)\|_0^2
    + 2\tau\langle D_\tau R^{\rm n.l.}_k,D_\tau e^*_k\rangle_H\\
    &\le c_a \tau\|D_\tau R_k\|_{H^*}^2
    + 6\tau^2\|\nabla D_\tau P(t_k)\|_0^2
    + \frac{\tau}{2c_a}(\|D_\tau e^*_k\|_1^2 + \|D_\tau e^*_{k-1}\|_1^2)\\
    &+ 2c_{\frac{1}{4c_a}}\tau\sum_{i=0}^1
    \left(\|e^*_{k-i-1}\|_1^2
    +\tau\left\|\frac{\partial u}{\partial t}\right\|_{L^2(t_{k-i-1},t_{k-i};\Hb^d)}^2\right).
\end{aligned}\]
Summing up for $k=2,3,\ldots,m$ with an arbitrary natural number 
$m \le N$, by Lemmas \ref{lem_R} and \ref{lem_diffdt}, it holds that 
\begin{align*}\begin{aligned}
    &\|D_\tau e_m\|_0^2 + \frac{\tau}{2c_a}\|D_\tau e^*_m\|_1^2
    + \tau\sum_{k=2}^m\frac{1}{2\tau}\|D_\tau e^*_k-D_\tau e_k\|_0^2\\
    \le\,& \|D_\tau e_1\|_0^2 + \frac{\tau}{2c_a}\|D_\tau e^*_1\|_1^2
    + \tau\sum_{k=2}^m(c_a\|D_\tau R_k\|_{H^*}^2
    + 6\tau\|D_\tau P(t_k)\|_1^2)\\
    &+ 4c_{\frac{1}{4c_a}}\tau\sum_{k=1}^m
    \left(\|e^*_k\|_1^2
    +\tau\left\|\frac{\partial u}{\partial t}\right\|_{L^2(t_{k-1},t_k;\Hb^d)}^2\right)\\
    \le\,& c_2\biggl\{\|D_\tau e_1\|_0^2
    + \tau\|D_\tau e^*_1\|_1^2
    + \tau^2\left\|\frac{\partial^3 u}{\partial t^3}\right\|_{L^2(H^*)}^2
    + \tau\left\|\frac{\partial P}{\partial t}\right\|_{L^2(H^1)}^2
    + \|\bar{e}^*_\tau\|_{L^2(H^1)}^2
    + \tau^2\left\|\frac{\partial u}{\partial t}\right\|_{L^2(H^1)}^2 \biggl\},
\end{aligned}\end{align*}
where $c_2:=\max\{2^{-1}c_a^{-1}, c_a, 6, 4c_{\frac{1}{4c_a}}\}$.
Hence, by Lemma \ref{lem_firstDt}, there exist two constants $c_3>0$ 
such that for all $0<\tau<\tau_1$,
\begin{align}\label{ineq_debdd}\begin{aligned}
    \max_{k=1,\ldots,N}\|D_\tau e_k\|_0^2
    + \tau\sum_{k=2}^N\frac{1}{\tau}\|D_\tau e_k-D_\tau e^*_k\|_0^2
    \le c_3.
\end{aligned}\end{align}

To use Lemma \ref{lem_AN} for (\ref{eq_dt}),
we set $(w_k,r_k)=T(D_\tau e_k)$ for all $k=1,2,\ldots,N$.
By Hypothesis \ref{hyp_stokes} and (\ref{ineq_debdd}), 
there exists a constant $c_4>0$ such that 
for all $0<\tau<\tau_1$ and $k=1,2,\ldots,N$,
$
    \|w_k\|_2\le c_4,
$
and hence, for all $0<\tau<\tau_1$, $\delta>0$ and $k=2,3,\ldots,N$,
\[\begin{aligned}
    -\tau d(\upre_{k-2},D_\tau e^*_k,w_k)
    &= \tau d(\upre_{k-2},w_k,D_\tau e^*_k)
    \le c_d\tau\|\upre_{k-2}\|_1\|w_k\|_2\|D_\tau e^*_k\|_0\\
    &\le c_dc_4 \tau \|\upre_{k-2}\|_1\|D_\tau e^*_k\|_0
    = \delta\|D_\tau e^*_k\|_0^2
    + \frac{c_d^2c_4^2}{4\delta} \tau^2\|\upre_{k-2}\|_1^2,\\  
    -\tau d(D_\tau u^*_{k-1},e^*_k,w_k)
    &= -\tau d(D_\tau u(t_{k-1}),e^*_k,w_k)
    - \tau d(D_\tau e^*_{k-1},e^*_k,w_k)\\
    &\le c_d\tau \|D_\tau u(t_{k-1})\|_1\|e^*_k\|_1\|w_k\|_1
    + c_d \|D_\tau e^*_{k-1}\|_0\|e^*_k\|_1\|w_k\|_2\\
    &\le c_d c_1\tau \|D_\tau u(t_{k-1})\|_1\|w_k\|_1
    + c_d c_4\|D_\tau e^*_{k-1}\|_0\|e^*_k\|_1\\
    &\le \frac{\delta}{2}\|w_k\|_1^2
    + \frac{c_d^2c_1^2}{2\delta}\tau^2\|D_\tau u(t_{k-1})\|_1^2
    + \delta\|D_\tau e^*_{k-1}\|_0^2
    + \frac{c_d^2c_4^2}{4\delta}\|e^*_{k-1}\|_1^2.    
\end{aligned}\]
By (\ref{eq_dtnl}) and (\ref{ineq_diffnl}) with $\varphi := w_k$, we have
\[\begin{aligned}
    \langle D_\tau R^{\rm n.l.}_k,w_k\rangle_H
    &\le \delta(\|D_\tau e^*_{k-1}\|_0^2+\|D_\tau e^*_k\|_0^2)\\
    &\quad + \tilde{c}_\delta\left\{\|w_k\|_1^2
    + \tau^2\|\upre_{k-2}\|_1^2
    + \sum_{i=0}^1
    \left(\|e^*_{k-i-1}\|_1^2
    + \tau^2\|D_\tau u(t_{k-i})\|_1^2\right)\right\},
\end{aligned}\]
where $\tilde{c}_\delta:=\max\{\delta,
    \delta^{-1}c_d^2(2c_{\rm max}^2 + 2^{-1}c_1^2 + 4^{-1}c_4^2)\}$.
By Lemmas \ref{lem_AN} and \ref{lem_diffdt}, 
there exist two constants $0 < \tau_2 \le \tau_1$ and $c_5>0$ such that
for all $0 < \tau < \tau_2$,
\[\begin{aligned}
    \tau\sum_{k=2}^N\|D_\tau e_k\|_0^2
    \le\,& c_5\biggl(\|D_\tau e_1\|_{V^*}^2
    + \tau\|D_\tau e^*_1\|_0^2
    + \tau^2\|\bar{u}^*_\tau\|_{L^2(H^1)}^2
    + \|\bar{e}^*_\tau\|_{L^2(H^1)}^2\\
    &+ \tau^2\left\|\frac{\partial u}{\partial t}\right\|_{L^2(H^1)}^2
    + \tau\sum_{k=2}^N(\|D_\tau e_k-D_\tau e^*_k\|_0^2
    + \|D_\tau R_k\|_{H^*}^2)\biggl).
\end{aligned}\]
Hence, by Theorems \ref{thm_stab}, \ref{thm_conv1}, 
Lemmas \ref{lem_R}, \ref{lem_firstDt}, and (\ref{ineq_debdd}), 
it holds that for all $0<\tau<\tau_2$,
\[
    \left\|\frac{\partial\hat{e}_\tau}{\partial t}\right\|_{L^2(L^2)}^2
    =\tau\|D_\tau e_1\|_0^2+\tau\sum^N_{k=2}\|D_\tau e_k\|_0^2
    \le c_7\tau
\]
for a constant $c_7>0$, where we have used 
$\frac{\partial\hat{e}_\tau}{\partial t}=(\overline{D_\tau e})_\tau$
on $(t_{k-1},t_k)$ for all $k=1,2,\ldots,N$
\qed
\medskip

Finally, we prove the second inequality of Theorem \ref{thm_conv2}.

\noindent{\it Proof of the second inequality of Theorem \ref{thm_conv2}.}
By (\ref{eq_error}) and Lemmas \ref{lem_infsup}, \ref{lem_a},
there exists a constant $c_1 > 0$ such that for all $k=1,2,\ldots,N$,
\begin{align*}\begin{aligned}
    \|q_k\|_0\!
    &\le\! c_1\sup_{0\ne\varphi\in H}\!\!
    \frac{|(q_k,\dvg\varphi)|}{\|\varphi\|_1}
    \!=\! c_1\sup_{0\ne\varphi\in H}\!\!\frac{\left|(D_\tau e_k,\varphi) 
    \!+\! a(e^*_k,\varphi)
    \!-\! \langle R_k \!+\! R^{\rm n.l.}_k, \varphi\rangle_H\right|}{\|\varphi\|_1}\\
    &\le c_1\left(\|D_\tau e_k\|_0 + c_a\|e^*_k\|_1 
    + \|R_k\|_{H^*} + \|R^{\rm n.l.}_k\|_{H^*}\right).
\end{aligned}\end{align*}
By Hypothesis \ref{hyp_reg_sol} and Theorem \ref{thm_conv1},
there exist two constants $\tau_1, c_2>0$ 
such that $\|u(t_k)\|_2$, $\tau^{-1/2}\|\bar{e}^*_\tau\|_{L^2(H^1)}\le c_2$
for all $0<\tau<\tau_1$ and $k=0,1,\ldots,N$.
By Lemma \ref{lem_diffdt}, it holds that 
for all $0<\tau<\tau_1$, $k=1,2,\ldots,N$ and $\varphi\in H$,
\[\begin{aligned}
    |\langle R^{\rm n.l.}_k,\varphi\rangle_H|
    &= |-d(e^*_{k-1},u(t_k),\varphi)
    - d(e^*_{k-1},e^*_k,\varphi)
    - d(u(t_{k-1}),e^*_k,\varphi)
    + \tau d(D_\tau u(t_k),u(t_k),\varphi)|\\
    &\le\! c_d\bigl(\|e^*_{k-1}\|_1\|u(t_k)\|_1
    \!+\! \|e^*_{k-1}\|_1\|e^*_k\|_1
    \!+\! \|u(t_{k-1})\|_1\|e^*_k\|_1
    \!+\! \tau \|D_\tau u(t_k)\|_0\|u(t_k)\|_2\bigr)\|\varphi\|_1\\
    &\le c_d c_2\left(\|e^*_{k-1}\|_1
    + 2\|e^*_k\|_1
    + \tau \|D_\tau u(t_k)\|_0\right)\|\varphi\|_1\\
    &\le c_d c_2\left(\|e^*_{k-1}\|_1
    + 2\|e^*_k\|_1
    + \sqrt{\tau} \left\|\frac{\partial u}{\partial t}\right\|_{L^2(t_{k-1},t_k,\Lo^d)}\right)\|\varphi\|_1,
\end{aligned}\]
where we have used $\|e_k\|_1\le c_2$ for all $k=0,1,\ldots,N$.
Hence, we have for all $0<\tau<\tau_1$ and $k=1,2,\ldots,N$,
\begin{align*}\begin{aligned}
    \|q_k\|_0
    \!\le\! c_3\left(\!\|D_\tau e_k\|_0 \!+\! \|e^*_{k-1}\|_1
    \!+\! \|e^*_k\|_1
    \!+\! \sqrt{\tau} \left\|\frac{\partial u}{\partial t}
        \right\|_{L^2(t_{k-1},t_k,\Lo^d)}
    \!\!+\! \|R_k\|_{H^*}\!\right)
\end{aligned}\end{align*}
for a constant $c_3>0$.
By Lemmas \ref{lem_R} and \ref{lem_dedt}, 
there exist three constants $\tau_2, c_4, c_5 > 0$
such that for all $0<\tau<\tau_2\le\tau_1$,
\[\begin{aligned}
    \|\bar{P}_\tau-P_\tau\|_{L^2(L^2)}^2
    &\le c_4\left(\left\|\frac{\partial \hat{e}_\tau}{\partial t}\right\|_{L^2(L^2)}^2
    + \|\bar{e}^*_\tau\|_{L^2(H^1)}^2
    + \tau^2\left\|\frac{\partial u}{\partial t}\right\|_{L^2(L^2)}^2
    + \|\bar{R}_\tau\|_{L^2(H^*)}^2\right)\\
    &\le c_4\left(\left\|\frac{\partial \hat{e}_\tau}{\partial t}\right\|_{L^2(L^2)}^2
    + c_2^2 \tau
    + \tau^2\left\|\frac{\partial u}{\partial t}\right\|_{L^2(L^2)}^2
    + \frac{\tau^2}{3}
    \left\|\frac{\partial^2 u}{\partial t^2}\right\|_{L^2(H^*)}^2\right)
    \le c_5\tau.
\end{aligned}\]
Therefore, by Lemma \ref{lem_diffdt}, we conclude the proof:
\[
    \|P-\bar{P}_\tau\|_{L^2(L^2)}
    \!\le\! \|P-P_\tau\|_{L^2(L^2)} + \|P_\tau-\bar{P}_\tau\|_{L^2(L^2)}
    \!\le\! \sqrt{\tau}\left(\left\|\frac{\partial P}{\partial t}\right\|_{L^2(L^2)}
    + \sqrt{c_5}\right).
\]
\qed

\section{Numerical examples}\label{sec_num}

\begin{figure}[htpb]
  \begin{minipage}{0.47\hsize}\centering
    \includegraphics[width=0.9\linewidth]{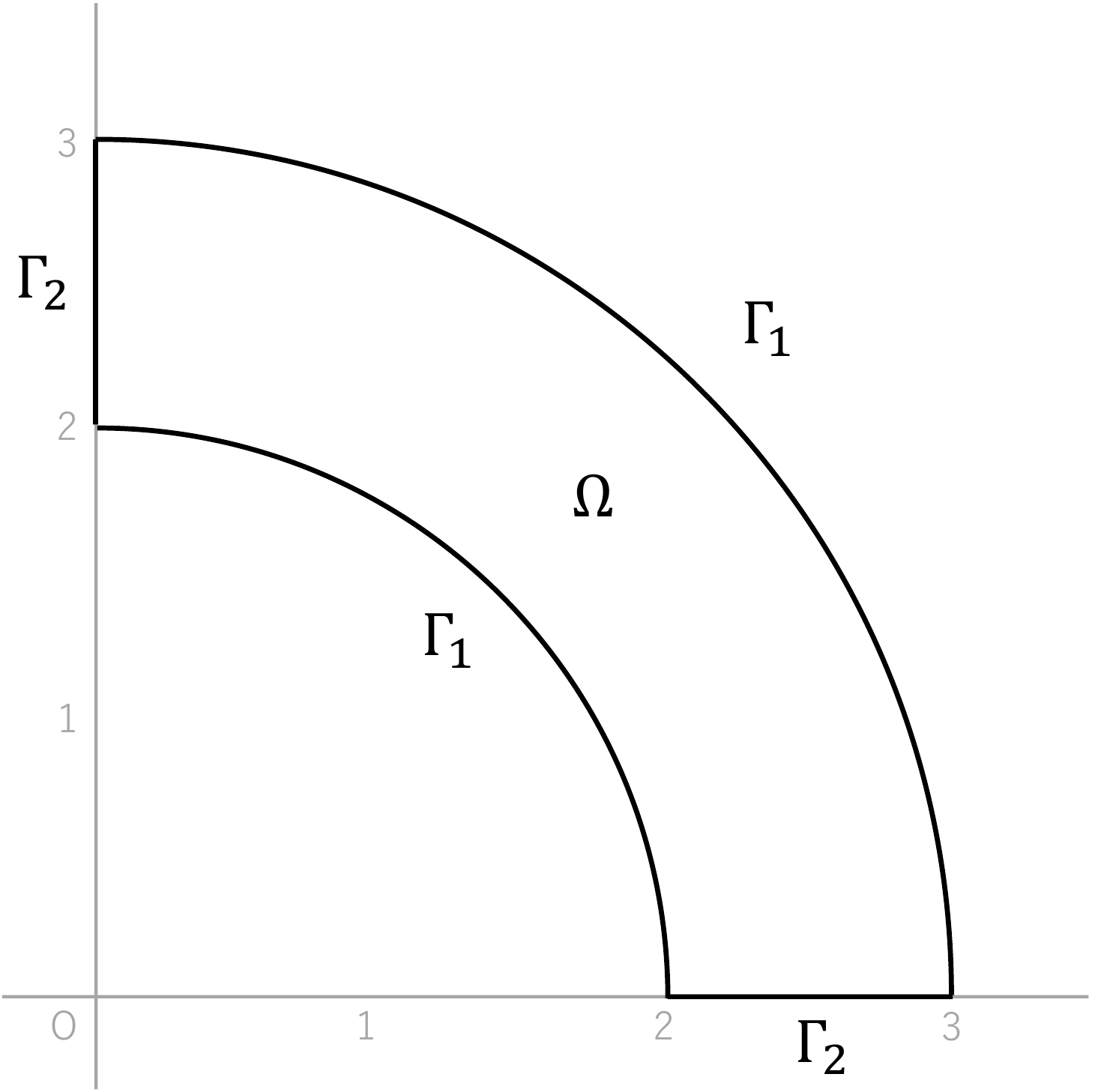}
      \end{minipage}
    \centering\begin{minipage}{0.47\hsize}\centering
    \includegraphics[width=0.9\linewidth]{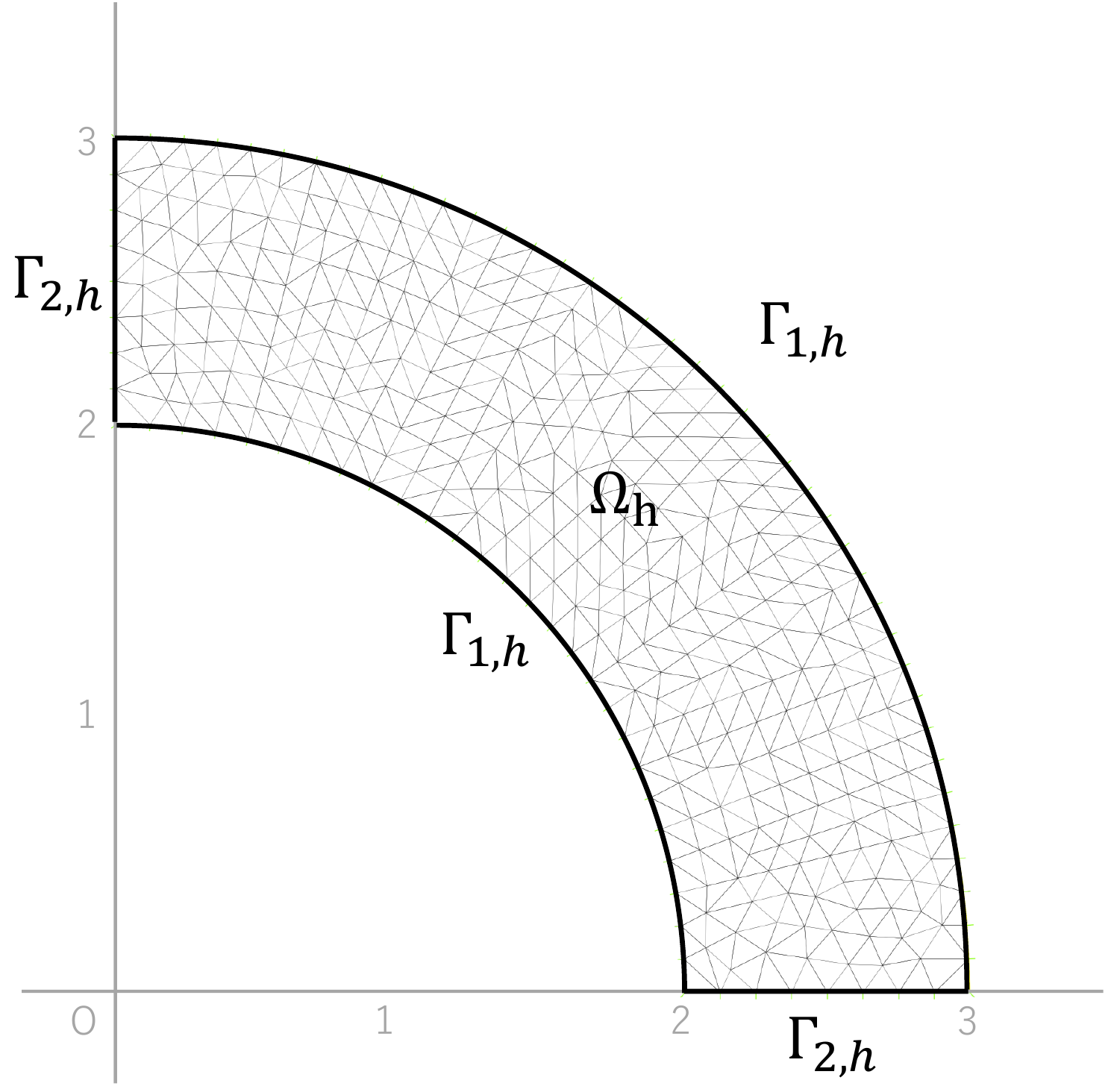}
    \end{minipage}
  \caption{The domain $\Omega$ with the boundary $\Gamma_1,\Gamma_2$ (left), 
  and $\Omega_h, \Gamma_{1,h}, \Gamma_{2,h}$ with mesh (right).}
  \label{fig_domain}       
\end{figure}

\begin{figure}[htpb]
  \begin{minipage}{0.47\hsize}\centering
    \includegraphics[width=0.9\linewidth, trim=30 0 20 0]{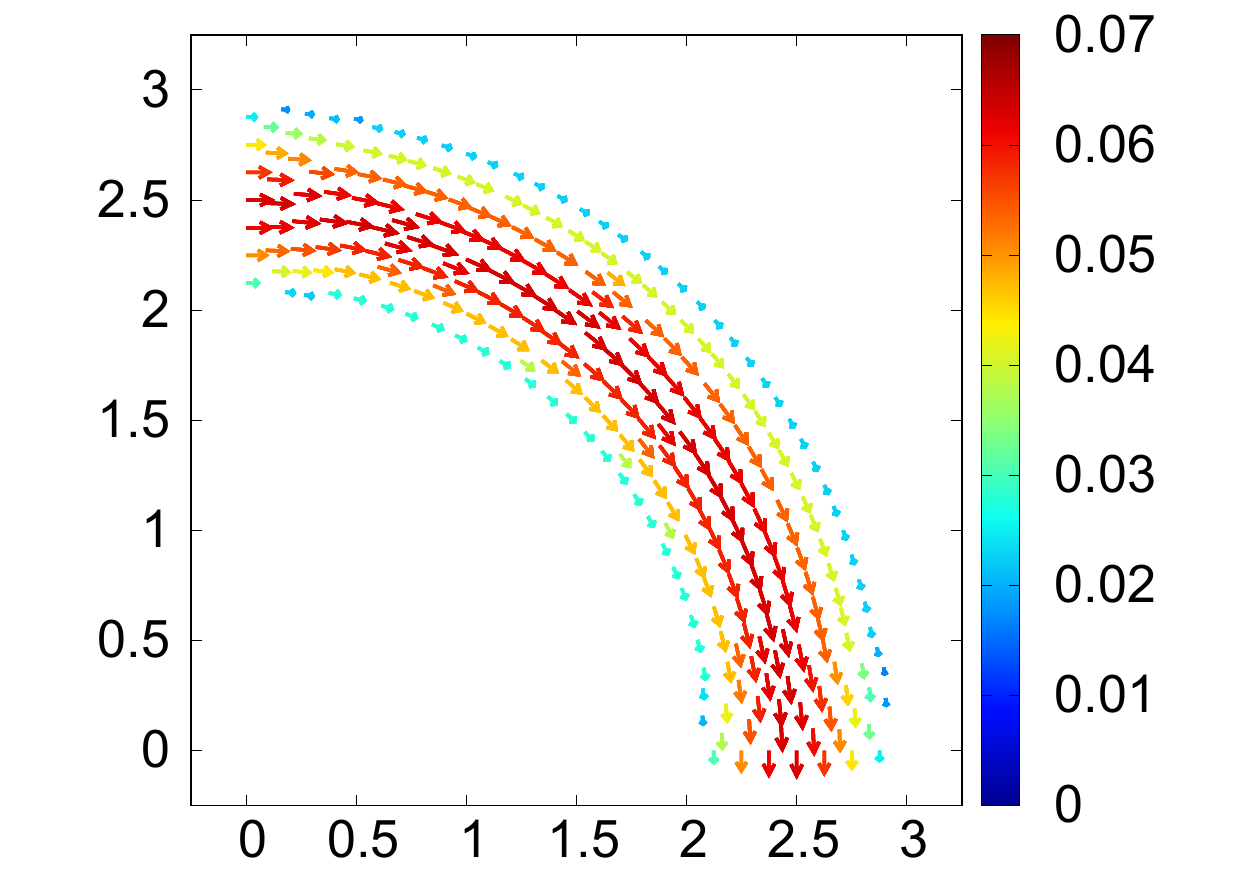}
      \end{minipage}
    \centering\begin{minipage}{0.47\hsize}\centering
    \includegraphics[width=0.9\linewidth, trim=20 0 30 0]{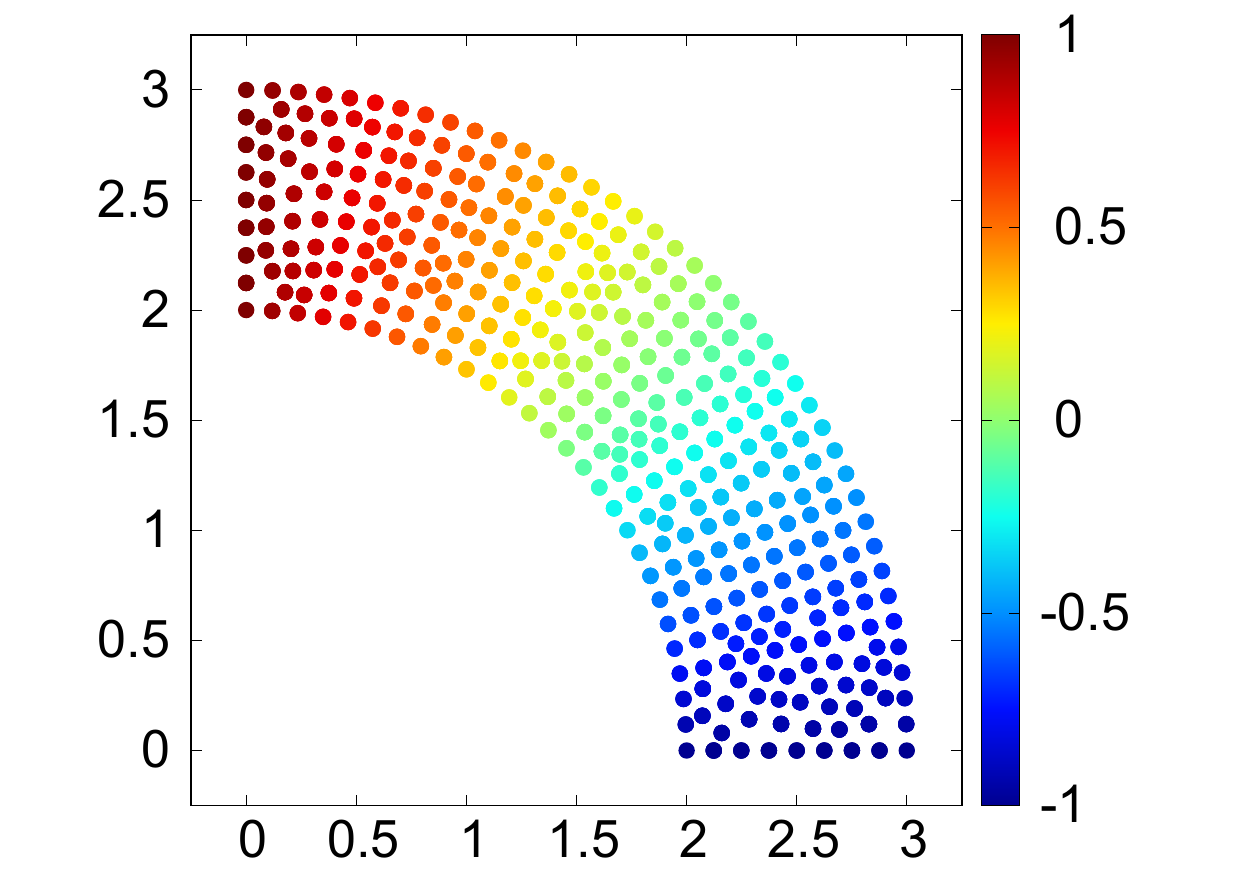}
    \end{minipage}
  \caption{The initial value $u_0$ of the velocity (left)
  and the pressure $p$ at $t=0$ (right).
  In the left figure, the color scale indicates 
  the length of $|u_0(\xi)|$ at each node $\xi$.}
  \label{fig_initial}       
\end{figure}

For our simulation, we set $T=1$ and 
\[\begin{aligned}
    \Omega &= \left\{(r\cos\theta, r\sin\theta) \in \R^2~\middle|~
        r_1<r<r_2, \theta_1<\theta<\theta_2\right\},\\
    \Gamma_1 &= \left\{(r\cos\theta, r\sin\theta) \in \R^2~\middle|~
    r\in \{r_1,r_2\}, \theta_1<\theta<\theta_2\right\},\\
    \Gamma_2 &= \left\{(r\cos\theta,r\sin\theta) \in \R^2~\middle|~
    r_1<r<r_2, \theta\in\{\theta_1,\theta_2\}\right\},
\end{aligned}\]
where $r_1:=2, r_2=3, \theta_1=0, \theta_2:=\pi/2$ (Fig. \ref{fig_domain}),
and define the following constants:
\[\begin{array}{c}
    {\DS p_{\rm in} := 1,\qquad 
    p_{\rm out} := -1,\qquad
    \alpha := \frac{p_{\rm in}-p_{\rm out}}{\theta_2-\theta_1},}\\
    {\DS C := \frac{1}{2}r_1^2 r_2^2\frac{\log\theta_2-\log\theta_1}{r_2^2-r_1^2},\quad
    D := -\frac{1}{2}\frac{r_2^2\log r_2-r_1^2\log r_1}{r_2^2-r_1^2}.}
\end{array}\]
The following functions 
\[
    u(x,y,t) := \left(\begin{array}{c}
        {\DS U(r)e^{-t}\sin\theta}\\
        {\DS -U(r)e^{-t}\cos\theta}
    \end{array}\right),\qquad
    p(x,y,t):= p_0(\theta)e^{-t},
\]
where $(r,\theta)=(r(x,y),\theta(x,y))$ are the polar coordinates and 
\[\begin{aligned}
    U(r) = \alpha\left(\frac{1}{2}r\log r + \frac{C}{r} + Dr\right),\qquad
    p_0(\theta) = \frac{
        p_{\rm in}(\theta-\theta_1) + p_{\rm out}(\theta_2-\theta)}{
        \theta_2-\theta_1},
\end{aligned}\]
satisfy (\ref{original}) with $\nu=\rho=1$ and 
\begin{align*}
    f(x,y,t)
    := \left(\begin{array}{c}
        {\DS -\frac{U^2(r)}{r} e^{-2t}\cos\theta
        -U(r) e^{-t}\sin\theta}\\[8pt]
        {\DS -\frac{U^2(r)}{r} e^{-2t}\sin\theta
        +U(r) e^{-t}\cos\theta}
    \end{array}\right)
    = \left\{\frac{\partial u}{\partial t}
    + (u\cdot\nabla)u\right\}(x,y,t),
\end{align*}
\begin{align*}
    p^b(x,y,t)
    := p_0(\theta)e^{-t} + \frac{U^2(r)}{2}e^{-2t},\qquad
    u_0(x,y) := \left(\begin{array}{c}
        {\DS U(r)\sin\theta}\\
        {\DS -U(r)\cos\theta}
    \end{array}\right).
\end{align*}
Fig. \ref{fig_initial} shows the initial value $u_0$ of the velocity
and the pressure $p$ at $t=0$.

We introduce a domain $\Omega_h$ to approximate the domain $\Omega$,
with boundary $\partial\Omega_h = \Gamma_{1,h} \cup \Gamma_{2,h}$
(Fig. \ref{fig_domain}).
We also introduce a regular triangulation $\mathcal{T}_h$ to $\Omega_h$,
with $h = \max_{K \in \mathcal{T}_h} \text{diam}(K)$
and $\overline{\Omega_h}=\cup_{K \in \mathcal{T}_h} \overline{K}$.
To consider the P2 and P1 element approximation 
for velocity and pressure, respectively,
we define the function spaces: for $i=1,2$,
\[\begin{aligned}
    X_h^i &:= \left\{ \psi_h \in C(\overline{\Omega_h}) ~\middle|~
        \varphi_h|_K \in P_i(K), \forall K \in \mathcal{T}_h
    \right\},\\
    H_h &:= \left\{ \varphi_h \in (X_h^2)^2 ~\middle|~
        \varphi_h = 0 \text{ on }\Gamma_{1,h},~
        \varphi_h\times n_h = 0 \text{ on }\Gamma_{2,h}
    \right\},\\
    Q_h &:= \left\{ \psi_h \in X_h^1 ~\middle|~
        \psi_h = 0 \text{ on }\Gamma_{2,h}
    \right\},
\end{aligned}\]
where $P_i(K)$ is the set of polynomials of degree $i$ or less on $K$
and $n_h$ is the unit outward normal vector for $\Gamma_{2,h}$.
Here, since $\Gamma_{2,h}$ is flat, 
the normal component of $\varphi_h\in H_h$ is not determined.
If $\Gamma_{2,h}$ is not flat, 
then $n_h$ is discontinuous on $\Gamma_{2,h}$
and $\varphi_h=0$ on $\Gamma_{2,h}$ (cf. \cite{BCPS17}).
Let $\Pi_h^i:C(\overline{\Omega_h}) \rightarrow X_h^i ~ (i=1,2)$
be the Lagrange interpolation operator (on each triangle).
By replacing $u_{k-1}$ in the first equation of (\ref{weak_FS})
with the third equation of (\ref{weak_FS}) at the previous step
(Remark \ref{rem_equiv}), we consider the following discrete problem:
\begin{Prob}\label{def_discrFS}
  For all $k=1,2,\ldots,N$, 
  find $(\upre_k,P_k)\in H_h\times X_h^1$
  such that $P_k - \Pi_h^1 p^b(t_k) \in Q_h$ and 
  for all $\varphi\in H_h$ and $\psi\in Q_h$,
  \begin{align}\label{weak_discrFS}\left\{\begin{array}{l}
      {\DS \frac{1}{\tau}(\upre_k-\upre_{k-1},\varphi)
      + a(\upre_k,\varphi)
      + d(\upre_{k-1},\upre_k,\varphi)
      + (\nabla P_{k-1},\varphi)}
      ={\DS (f(t_k),\varphi),}\\[8pt]
      {\DS \tau(\nabla P_k,\nabla\psi)}
      ={\DS -(\dvg\upre_k,\psi),}\\[4pt]
  \end{array}\right.\end{align}
  where $P_0:=0$.
\end{Prob}
For all $k=1,2,\ldots,N$, we set $u_k := \upre_k - \tau\nabla P_k$.
See \cite{GQ97,GQ98} for the details on $u_k$ and its divergence.

\begin{figure}[htpb]
  \begin{minipage}{0.47\hsize}\centering
    \includegraphics[width=\linewidth, trim=100 0 80 0]{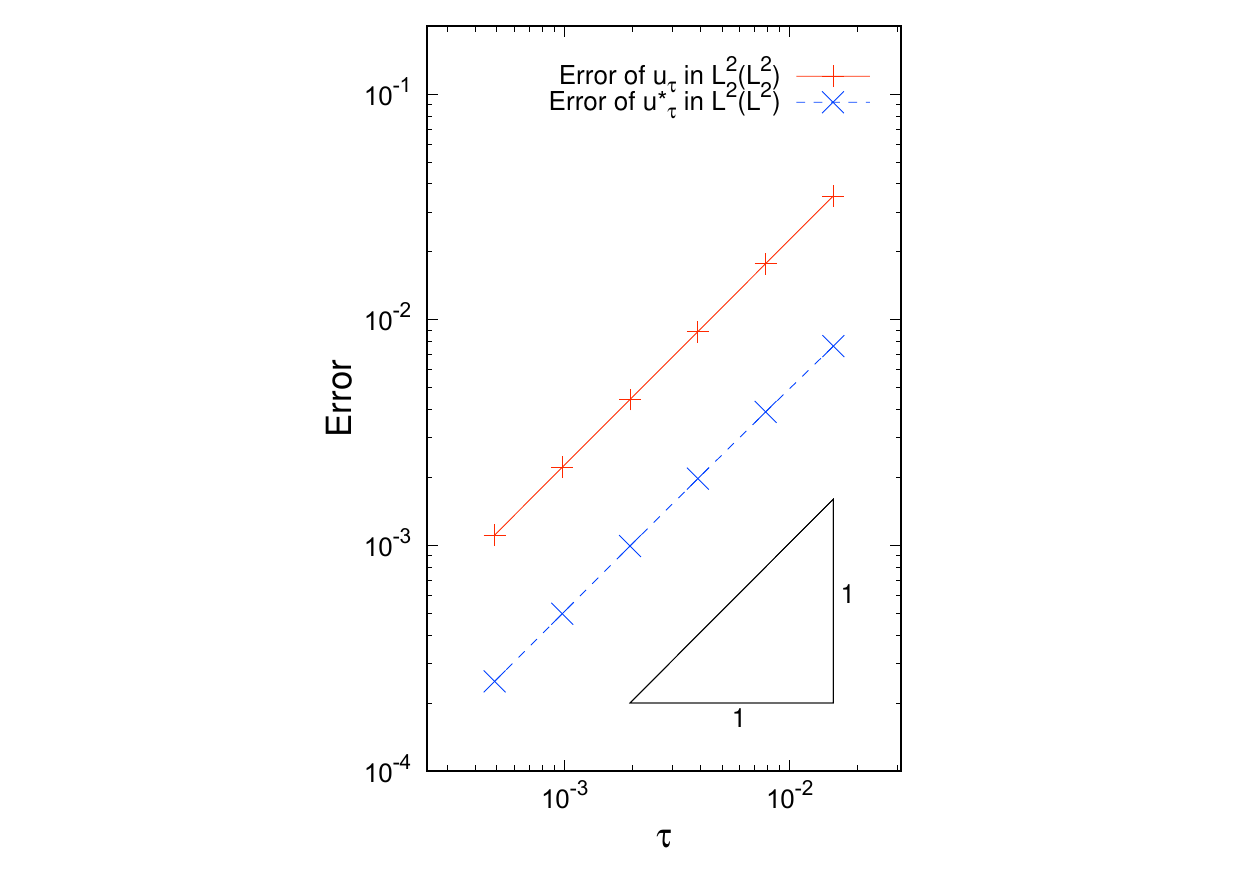}
      \end{minipage}
    \centering\begin{minipage}{0.47\hsize}\centering
    \includegraphics[width=\linewidth, trim=80 0 100 0]{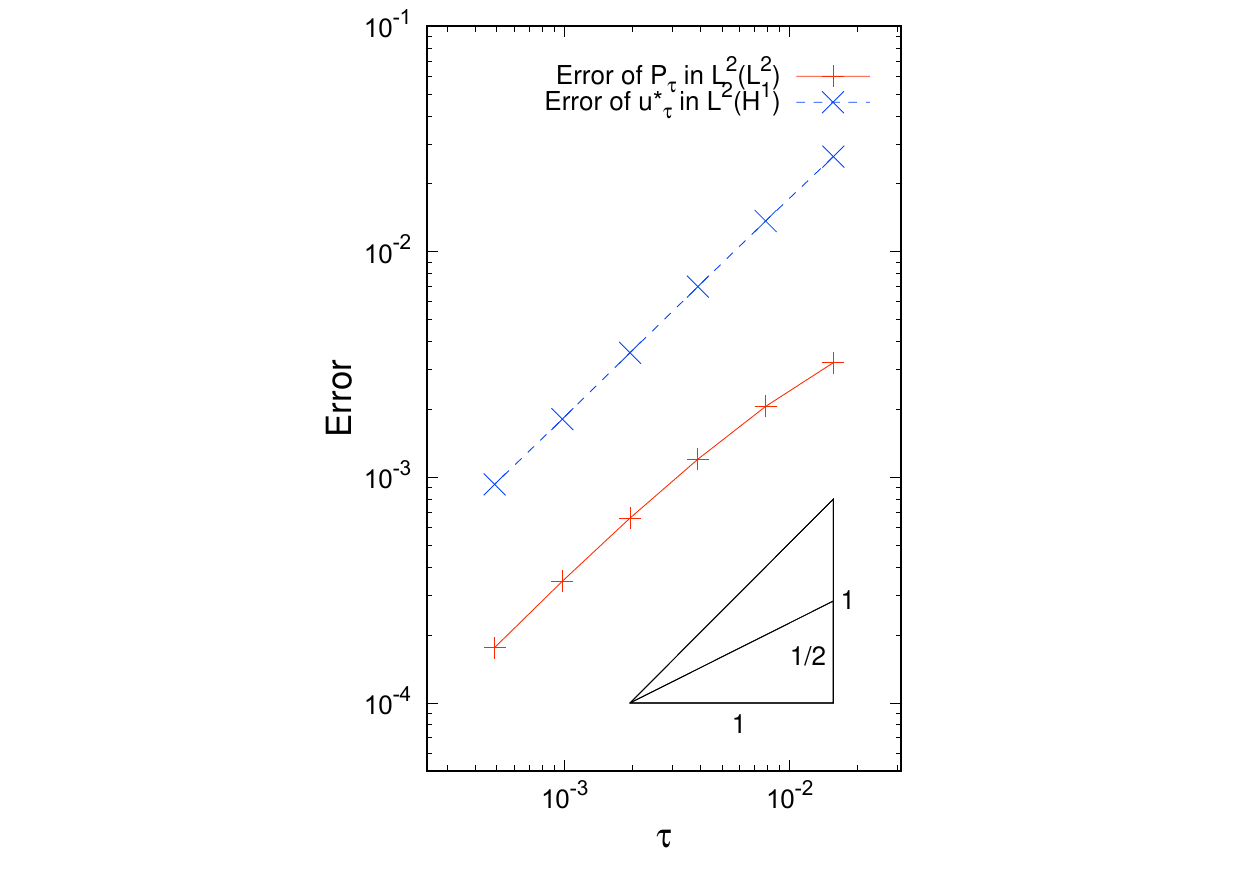}
    \end{minipage}
  \caption{The errors of logscale:
  $\|\bar{u}_\tau - \Pi_h^2 u_\tau\|_{L^2(L^2(\Omega_h)^d)}$,
  $\|\bar{u}^*_\tau - \Pi_h^2 u_\tau\|_{L^2(L^2(\Omega_h)^d)}$ (left),
  $\|\bar{P}_\tau - \Pi_h^1 P_\tau\|_{L^2(L^2(\Omega_h))}$,
  and $\|\bar{u}^*_\tau - \Pi_h^2 u_\tau\|_{L^2(H^1(\Omega_h)^d)}$ (right).
  The triangles show the slope of $O(\tau)$ and $O(\sqrt{\tau})$.}
  \label{fig_error}       
\end{figure}

On a mesh with $h = 2^{-6}$,
we solve the problems (\ref{weak_discrFS}) numerically by
using the software FreeFEM \cite{FreeFem}.
We compute the error estimates between 
the numerical solutions of (\ref{weak_discrFS})
and the interpolation $(\Pi_h^2 u, \Pi_h^1P)$ 
of the exact solution $(u,P)$, where $P:=p+|u|^2/2$.
In Fig. \ref{fig_error}, the numerical errors 
$\|\bar{u}_\tau - \Pi_h^2 u_\tau\|_{L^2(L^2(\Omega_h)^d)}$,
$\|\bar{u}^*_\tau - \Pi_h^2 u_\tau\|_{L^2(L^2(\Omega_h)^d)}$,
$\|\bar{P}_\tau - \Pi_h^1 P_\tau\|_{L^2(L^2(\Omega_h))}$,
and $\|\bar{u}^*_\tau - \Pi_h^2 u_\tau\|_{L^2(H^1(\Omega_h)^d)}$
are presented.
One can observe that 
$\|\bar{u}_\tau - \Pi_h^2 u_\tau\|_{L^2(L^2(\Omega_h)^d)}$ and
$\|\bar{u}^*_\tau - \Pi_h^2 u_\tau\|_{L^2(L^2(\Omega_h)^d)}$
are almost of first order in $\tau$ and that 
$\|\bar{P}_\tau - \Pi_h^1 P_\tau\|_{L^2(L^2(\Omega_h))}$
is of 0.5th order in $\tau$,
as expected from Theorem \ref{thm_conv2}.
Furthermore, the error
$\|\bar{u}^*_\tau - \Pi_h^2 u_\tau\|_{L^2(H^1(\Omega_h)^d)}$
is almost of first order in $\tau$,
which is better than the theoretically predicted rate 
(Theorem \ref{thm_conv1}).

\section{Conclusion}\label{sec_conclusion}
We have proposed a new projection method 
for Navier--Stokes equations (\ref{original}) 
with a total pressure boundary condition. 
We have shown the stability of the projection method 
in Theorem \ref{thm_stab} and established 
error estimates for the velocity and the pressure 
in suitable norms between the solution 
to (\ref{weak_original}) and (\ref{weak_FS}) 
in Theorems \ref{thm_conv1} and \ref{thm_conv2}. 
The convergence rates are the same as the case of 
the usual full-Dirichlet boundary condition for velocity \cite{Shen92,Shen94}. 
The traction boundary condition is often used to apply 
Dirichlet boundary conditions for pressure; however, 
the convergence rates are worse than our case
(Compare \cite{GMS05} and \cite{GS04}).

The projection method is still evolving, 
and many high-convergence methods have been proposed \cite{GMS06}. 
The application of the boundary conditions proposed in this paper 
to these methods will be a focus of our future works. 
As another future direction, 
the case that $\Gamma_2$ is not flat in numerical calculations 
is an important problem (cf. \cite{BCPS17}). 
In addition, since the nonlinear term $(\nabla\times u)\times u$ 
is different from the standard advection term $(u\cdot\nabla)u$, 
it cannot be applied to methods using the Lagrangian coordinates, 
such as the characteristic curve method and particle methods;
this problem remains open for further study. 

\section*{Acknowledgments}

This work was supported by JSPS KAKENHI Grant Number 19J20514.

\bibliographystyle{spmpsci}      
\bibliography{PSwPressure}   

\appendix
\section*{Proofs of Lemmas \ref{lem_stokes}, \ref{lem_nonlin},
\ref{lem_R}, \ref{lem_diffdt}, and Corollary \ref{cor_existsol}}

The purpose of this appendix is to provide 
the proofs of Lemmas \ref{lem_stokes}, \ref{lem_nonlin},
\ref{lem_R}, \ref{lem_diffdt}, and Corollary \ref{cor_existsol}.
The continuity of the operator $T$ follows from Lemmas \ref{Stokes},
\ref{lem_infsup}, and \ref{lem_a}. By using Lemma \ref{lem_a} again,
we prove the second inequality of Lemma \ref{lem_stokes}.

\setcounter{equation}{0}
\renewcommand{\theequation}{A.\arabic{equation}}
\setcounter{figure}{0}
\renewcommand{\thefigure}{\Alph{section}.\arabic{figure}}
\setcounter{table}{0}
\renewcommand{\thetable}{\Alph{section}.\arabic{table}}

\medskip
\noindent{\it Proof of Lemma \ref{lem_stokes}.}
By Lemmas \ref{Stokes}, \ref{lem_infsup} and \ref{lem_a},
there exists a unique solution $(w,r)\in H\times\Lo$ to (\ref{eq_stokes}) 
for all $e\in\Lo^d$ and $T$ is a continuous operator;
\[
    \|w\|_1+\|r\|_0\le c_1\|e\|_{H^*} \le c_1\|e\|_0
\]
for a constant $c_1>0$ independent of $e$.
It is easy to check that $T$ is a linear operator.

Next, we show the second inequality of Lemma \ref{lem_stokes}.
By the first equation of (\ref{eq_stokes}), 
it holds that for all $\varphi\in V$,
\[
    a(w,\varphi)=(r,\dvg\varphi)+(e,\varphi)=(e,\varphi).
\]
By Lemma \ref{lem_a}, we have
\[
    \frac{1}{c_a}\|w\|_1^2
    \le\|w\|_a^2=a(w,w)=(e,w)\le\|e\|_{V^*}\|w\|_1,
\]
which implies that $\|w\|_1\le c_a\|e\|_{V^*}$.
On the other hand, by Lemma \ref{lem_a},
it holds that for all $e\in\Lo^d$,
\[
    \|e\|_{V^*}
    =\sup_{0\ne\varphi\in V}\frac{(e,\varphi)}{\|\varphi\|_1}
    =\sup_{0\ne\varphi\in V}\frac{a(w,\varphi)}{\|\varphi\|_1}
    \le \sup_{0\ne\varphi\in V}\frac{c_a\|w\|_1\|\varphi\|_1}{\|\varphi\|_1}
    =c_a\|w\|_1.
\]
\qed

In order to prove Lemma \ref{lem_nonlin} and Corollary \ref{cor_existsol},
we define $\tilde{p}_d, \tilde{q}_d$ as
\[
    \tilde{p}_d := \frac{2p_d}{p_d+2}
    = \frac{1}{\frac{1}{2} + \frac{1}{p_d}},\qquad
    \tilde{q}_d := \frac{\tilde{p}_d}{1 - \tilde{p}_d}
    = \frac{1}{1 - \frac{1}{2} - \frac{1}{p_d}}.
\]
Here, since $p_2=2+\eps$ and $p_3=3$,
we have $(\tilde{p}_2, \tilde{q}_2) = (\frac{4+2\eps}{4+\eps}, 2+\frac{4}{\eps})$
and $(\tilde{p}_2, \tilde{q}_2) = (\frac{6}{5}, 6)$.
By the Sobolev embeddings \cite[Theorem III.2.33]{BF13},
it holds that $\Hb \subset L^{\tilde{q}_d}(\Omega)$,
$\Hb \subset L^{p_d}(\Omega)$,
$H^2(\Omega) \subset L^\infty(\Omega)$
and the embeddings are continuous.

\medskip
\noindent{\it Proof of Lemma \ref{lem_nonlin}.}

(i) For all $u \in L^{p_d}(\Omega)^d, v, w \in H$, we have
\begin{align*}
    |d(u,v,w)|
    &\le \into \left|u\cdot((w\cdot\nabla)v-(v\cdot\nabla)w
    +v\,\dvg w-w\,\dvg v)\right|dx\\
    &\le c_1\|u\|_{L^{p_d}}
    (\|w\|_{L^{\tilde{q}_d}}\|\nabla v\|_0
    +\|v\|_{L^{\tilde{q}_d}}\|\nabla w\|_0
    +\|v\|_{L^{\tilde{q}_d}}\|\dvg w\|_0
    +\|w\|_{L^{\tilde{q}_d}}\|\dvg v\|_0)\\
    &\le \tilde{c}_1\|u\|_{L^{p_d}}\|v\|_1\|w\|_1
\end{align*}
for two constants $c_1,\tilde{c}_1>0$,
which implies the third inequality of Lemma \ref{lem_nonlin}.

(ii) For all $u \in L^{p_d}(\Omega)^d, v \in H, w \in H \cap H^2(\Omega)^d$, 
we have
\begin{align*}
    |d(u,v,w)|
    &\le c_2\|u\|_0
    (\|w\|_{L^\infty}\|\nabla v\|_0
    +\|v\|_{L^{p_d}}\|\nabla w\|_{L^{\tilde{q}_d}}
    +\|v\|_{L^{p_d}}\|\dvg w\|_{L^{\tilde{q}_d}}
    +\|w\|_{L^\infty}\|\dvg v\|_0)\\
    &\le \tilde{c}_2\|u\|_0\|v\|_1\|w\|_2
\end{align*}
for two constants $c_2,\tilde{c}_2>0$.

(iii) For all $u \in \Hb^d, v \in H \cap H^2(\Omega)^d, w \in H$, 
we have
\[
    |d(u,v,w)|
    \le \into |((\nabla\times u)\times v)\cdot w|dx
    \le c_3\|\nabla\times u\|_0\|v\|_{L^\infty}\|w\|_0
    \le \tilde{c}_3\|u\|_1\|v\|_2\|w\|_0
\]
for two constants $c_3,\tilde{c}_3>0$.

(ix) For all $u \in H^2(\Omega)^d, v, w\in H$, we have 
\[
    |d(u,v,w)|
    \le c_4\|\nabla\times u\|_{L^{p_d}}\|v\|_{L^{\tilde{q}_d}}\|w\|_0
    \le \tilde{c}_4\|u\|_2\|v\|_1\|w\|_0
\]
for two constants $c_4,\tilde{c}_4>0$.
\qed

\medskip

Next, we prove Lemmas \ref{lem_R} and \ref{lem_diffdt}.

\medskip
\noindent{\it Proof of Lemma \ref{lem_R}.}
It holds that for all $\varphi\in H$ and $k=1,2,\ldots,N$, 
\[\begin{aligned}
  \langle R_k,\varphi\rangle_H
  &= \tau\int^1_0\left\langle
  s\frac{\partial^2 u}{\partial t^2}(t_{k-1}+s\tau),\varphi\right\rangle_H ds
  \le \tau\int^1_0
  s\left\|\frac{\partial^2 u}{\partial t^2}(t_{k-1}+s\tau)\right\|_{H^*}
  \|\varphi\|_1ds\\
  &\le
  \sqrt{\frac{\tau}{3}}\|\varphi\|_1
  \left\|\frac{\partial^2 u}{\partial t^2}\right\|_{L^2(t_{k-1},t_k;H^*)},
\end{aligned}\]
which implies that
\[
  \|\bar{R}_\tau\|_{L^2(H^*)}^2
  =\sum^N_{k=1}\tau\left(\sup_{0\ne\varphi\in H}
  \frac{\langle R_k,\varphi\rangle_H}{\|\varphi\|_1}\right)^2
  \le \frac{1}{3}\tau^2
  \left\|\frac{\partial^2 u}{\partial t^2}\right\|_{L^2(H^*)}^2.
\]

Next, we show the second inequality of the conclusion.
For all $\varphi\in H$ and $k=2,3,\ldots,N$, we have
\[\begin{aligned}
  \langle D_\tau R_k,\varphi\rangle_H
  &= \tau\int^1_0\int^1_0\left\langle
  s_1\frac{\partial^3 u}{\partial t^3}(t_{k-2}+s_1\tau+s_2\tau),
  \varphi\right\rangle_H ds_1 ds_2\\
  &\le \tau\int^1_0\int^1_0 s_1
  \left\|\frac{\partial^3 u}{\partial t^3}(t_{k-2}+s_1\tau+s_2\tau)\right\|_{H^*}
  \|\varphi\|_1 ds_1 ds_2\\
  &\le \tau\|\varphi\|_1\sqrt{\int^1_0\int^1_0 s_1^2 ds_1 ds_2}
  \sqrt{\int^1_0\int^1_0 \left\|\frac{\partial^3 u}{\partial t^3}
      (t_{k-2}+s_1\tau+s_2\tau)\right\|_{H^*}^2 ds_1 ds_2}\\
  &\le \tau\|\varphi\|_1\sqrt{\frac{1}{3}}
  \sqrt{\int^1_{-1}\int^2_0
      \left\|\frac{\partial^3 u}{\partial t^3}
      (t_{k-2}+\tilde{s}_1\tau)\right\|_{H^*}^2 \frac{1}{2} d\tilde{s}_1 d\tilde{s}_2}\\
  &= \sqrt{\frac{\tau}{3}}\|\varphi\|_1
  \left\|\frac{\partial^3 u}{\partial t^3}\right\|_{L^2(t_{k-2},t_k;H^*)},
\end{aligned}\]
where we have used the coordinate transformation 
$(s_1,s_2)\mapsto(\tilde{s}_1,\tilde{s}_2):=(s_1+s_2,-s_1+s_2)$.
Therefore, we obtain
\[
  \sum^N_{k=2}\tau\|D_\tau R_k\|_{H^*}^2
  \le \sum^N_{k=2}\tau\frac{\tau}{3}
  \left\|\frac{\partial^3 u}{\partial t^3}\right\|_{L^2(t_{k-2},t_k;H^*)}^2
  \le \frac{2}{3}\tau^2
  \left\|\frac{\partial^3 u}{\partial t^3}\right\|_{L^2(H^*)}^2.
\]
\qed

\medskip
\noindent{\it Proof of Lemma \ref{lem_diffdt}.}
It holds that for all $k = 1,2,\ldots,N$ and $t \in [t_{k-1}, t_k]$,
\[
    \|x(t_k) - x(t)\|_E
    \le \int^{t_k}_{t}
    \left\|\frac{\partial x}{\partial t}(s)\right\|_E ds
    \le \sqrt{t_k - t}
    \left\|\frac{\partial x}{\partial t}\right\|_{L^2(t_{k-1},t_k;E)},
\]
which implies that 
$\|x - x_\tau\|_{L^\infty(E)}
\le \sqrt{\tau}
\left\|\frac{\partial x}{\partial t}\right\|_{L^2(E)}$ and 
\[
    \|D_\tau x(t_k)\|_E
    =\frac{1}{\tau} \|x(t_k) - x(t_{k-1})\|_E
    \le \frac{1}{\sqrt{\tau}}
    \left\|\frac{\partial x}{\partial t}\right\|_{L^2(t_{k-1},t_k;E)}.
\]
On the other hand, we have
\[\begin{aligned}
    \|x-x_\tau\|_{L^2(E)}^2
    &=\sum^N_{k=1}\int^{t_k}_{t_{k-1}}\|x(t)-x(t_k)\|_E^2 dt\\
    &\le \sum^N_{k=1}\int^{t_k}_{t_{k-1}}(t_k-t)dt
        \left\|\frac{\partial x}{\partial t}\right\|_{L^2(t_{k-1},t_k;E)}^2
    = \frac{1}{2}\tau^2\left\|\frac{\partial x}{\partial t}\right\|_{L^2(E)}^2.
\end{aligned}\]
\qed

\medskip
We prove Corollary \ref{cor_existsol}
by using the boundedness from Theorem \ref{thm_stab}
and the Aubin--Lions compactness lemma.

\medskip
\noindent{\it Proof of Corollary \ref{cor_existsol}.}
By the first and third equations of (\ref{weak_FS}),
it holds that for all $v\in V$ and $k=1,2,\ldots,N$,
\[\begin{aligned}
    (D_\tau u_k,v)+a(\upre_k,v)+(g_k,v)+(h_k,\nabla v)
    =\langle f_k,v\rangle_H-(\nabla P_k,v)
    =\langle f_k,v\rangle_H-\int_{\Gamma_2}p^b_kv\cdot nds,
\end{aligned}\]
where $g_k$ and $h_k$ are defined\footnote{Here, it holds that 
for all $i,j=1,\ldots,d$ and $k=1,2,\ldots,N$,
\[
    (g_k)_i
    :=\sum_{l=1}^d\frac{\partial (\upre_k)_l}{\partial x_i}(\upre_{k-1})_l
    -(\upre_{k-1})_i\dvg\upre_k,\qquad
    (h_k)_{ij} := -(\upre_k)_i(\upre_{k-1})_j.
\]} by 
\[
    g_k := (\nabla\upre_k)^T \upre_{k-1}
    - \upre_{k-1}\dvg\upre_k,\qquad
    h_k := -\upre_k(\upre_{k-1})^T,
\]
which implies that for all $v\in V$ and $\theta\in C^\infty_0(0,T)$,
\begin{align}\label{eq_ghver}\begin{aligned}
    \int^T_0\left(\left(\frac{\partial \hat{u}_\tau}{\partial t},v\right)
    +a(\bar{u}^*_\tau,v)
    +(\bar{g}_\tau,v)+(\bar{h}_\tau,\nabla v)\right)\theta dt
    \!=\!\int^T_0\left(\langle \bar{f}_\tau,v\rangle_H
    - \int_{\Gamma_2}\bar{p}^b_\tau v\cdot nds\right)\theta dt.
\end{aligned}\end{align}
Here, $\bar{f}_\tau\rightarrow f$ strongly in $L^2(H^*)$
and $\bar{p}^b_\tau\rightarrow p^b$ strongly in $L^2(\Hb)$ as $\tau\rightarrow 0$.
By Theorem \ref{thm_stab} and Lemma \ref{lem_gen}, 
there exists a constant $c_1>0$ such that 
\begin{align}\label{ineq_bdd}\begin{aligned}
    \|\bar{u}_\tau\|_{L^\infty(L^2)}^2
    \!+\! \|\bar{u}^*_\tau\|_{L^\infty(L^2)}^2
    \!+\! \|\bar{u}^*_k\|_{L^2(H^1)}^2
    \!+\! \tau\left\|\frac{\partial\hat{u}_\tau}{\partial t}\right\|_{L^2(L^2)}^2
    \!\!+\! \frac{1}{\tau}\|\bar{u}_\tau-\bar{u}^*_\tau\|_{L^2(L^2)}^2
    \!\le\! c_1.
\end{aligned}\end{align}
In particular, it holds that 
\begin{align}\label{ineq_U0}
    \|\upre_1\|_0^2 + \tau\|\upre_1\|_1^2 
    + \|u_1-u_0\|_0^2 + \|u_1-\upre_1\|_0^2 
    \le c_1,
\end{align}
which implies that $\|\upre_1-u_0\|_0 
\le \|\upre_1-u_1\|_0 +  \|u_1-u_0\|_0 \le 2\sqrt{c_1}$.
Furthermore, by the first equation of (\ref{weak_FS})
and Lemmas \ref{lem_a}, \ref{lem_nonlin}, we have
\begin{align}\label{ineq_Ustar}\begin{aligned}
    \|\upre_1-u_0\|_{V^*}
    &= \sup_{0\ne v\in H}\frac{\tau}{\|v\|_1}
        \left|-a(\upre_1,v) - d(u_0,\upre_1,v)
        + \langle f_1, v\rangle_H\right|\\
    &\le c_a\tau\|\upre_1\|_1 + c_d\tau\|u_0\|_{L^{p_d}}\|\upre\|_1
    + \tau\|f_1\|_{H^*}
    \le c_2\sqrt{\tau}.
\end{aligned}\end{align}
where $c_2:=\sqrt{c_1}(c_a + c_d\|u_0\|_{L^{p_d}}) + \|f\|_{L^2(H^*)}$.
Let $u^\circ_0:=\upre_1$, $u^\circ_k:=\upre_k$ for all $k=1,2,\ldots,N$
and let $\hat{u}^\circ_\tau$ be 
the piecewise linear interpolant of $(u^\circ_k)_{k=0}^N \subset H$. 

From the uniform estimates (\ref{ineq_bdd}),
one can show that there exist a sequence $(\tau_k)_{k\in\N}$
and three functions $u\in L^2(H)\cap L^\infty(\Lo^d)\cap W^{1,4/p_d}(V^*)$
(in particular, $u\in C([0,T];V^*)$), 
$g\in L^{4/p_d}(L^{\tilde{p}_d}(\Omega)^d)$ and 
$h\in L^{4/p_d}(\Lo^{d\times d})$ such that 
$\tau_k\rightarrow 0$ and 
\begin{align}\label{convsbwh}
    \bar{u}^*_{\tau_k}\rightarrow u
    &\quad\text{weakly in }L^2(H),\\
    &\quad\begin{aligned}\label{convsbsl}
        \text{strongly in }L^2(\Lo^d),
    \end{aligned}\\
    \hat{u}^\circ_{\tau_k}\rightarrow u
    &\quad\begin{aligned}\label{convshsl}
        \text{strongly in }L^2(\Lo^d),
    \end{aligned}\\
    &\quad\begin{aligned}\label{convshsc}
        \text{strongly in }C([0,T];V^*),
    \end{aligned}\\
    \hat{u}_{\tau_k} \rightarrow u
    &\quad\begin{aligned}\label{convnhsl}
        \text{strongly in }L^2(\Lo^d),
    \end{aligned}\\
    &\quad\begin{aligned}\label{convnhwvs}
        \text{weakly in }W^{1,4/p_d}(V^*),
    \end{aligned}\\
    \bar{g}_{\tau_k} \rightharpoonup g
    &\quad\begin{aligned}\label{convgw}
        \text{weakly in }L^{4/p_d}(L^{\tilde{p}_d}(\Omega)^d),
    \end{aligned}\\
    \bar{h}_{\tau_k} \rightharpoonup h
    &\quad\begin{aligned}\label{convhw}
        \text{weakly in }L^{4/p_d}(\Lo^{d\times d}),
    \end{aligned}
\end{align}
as $k\rightarrow\infty$.
Here, we note that $\bar{u}^*_{\tau_k}$, $\hat{u}^\circ_{\tau_k}$
and $\hat{u}_{\tau_k}$ possess a common limit function.
Indeed, the weak convergence (\ref{convsbwh}) of $\bar{u}^*_\tau$ 
immediately follows from the uniform estimates (\ref{ineq_bdd}). 
Since we have $1/\tilde{p}_d = 1/2 + 1/p_d$,
$p_d/4 = 1/2 + p_d/(2\tilde{q}_d)$, and 
\[
    \|\bar{u}^*_\tau\|_{L^{2\tilde{q}_d/p_d}(L^{p_d})}
    \le \|\bar{u}^*_\tau\|_{L^2(L^{\tilde{q}_d})}^{p_d/\tilde{q}_d}
    \|\bar{u}^*_\tau\|_{L^\infty(L^2)}^{1-p_d/\tilde{q}_d}
    \le c_3\|\bar{u}^*_\tau\|_{L^2(H^1)}^{p_d/\tilde{q}_d}
    \|\bar{u}^*_\tau\|_{L^\infty(L^2)}^{1-p_d/\tilde{q}_d}
\]
for a constant $c_3>0$
(cf. \cite[Theorem II.5.5]{BF13})\footnote{
    Since $p_2=2+\eps$ and $p_3=3$, 
    we have $p_2/\tilde{q}_2 = \eps/2$ and $p_3/\tilde{q}_3=1/2$.
}, it holds that 
\[\begin{aligned}
    &\quad\|\bar{g}_\tau\|_{L^{4/p_d}(L^{\tilde{p}_d})}
    \le\left\{\tau\sum^N_{k=1}
    \left(\|\nabla\upre_k\|_{L^2}\|\upre_{k-1}\|_{L^{p_d}}
    +\|\upre_{k-1}\|_{L^{p_d}}\|\dvg\upre_k\|_{L^2}\right)^{4/{p_d}}\right\}^{{p_d}/4}\\
    &\le c_4\left(\tau\sum^N_{k=1}
    \|\upre_k\|_1^{4/{p_d}}\|\upre_{k-1}\|_{L^{p_d}}^{4/{p_d}}\right)^{{p_d}/4}
    \le c_4\left(\tau\sum^N_{k=1}\|\upre_k\|_1^2\right)^{1/2}
    \left(\tau\sum^N_{k=1}\|\upre_{k-1}\|_{L^{p_d}}^{2\tilde{q}_d/p_d}\right)^{p_d/(2\tilde{q}_d)}\\    
    &\le c_4
    \left(\|\bar{u}^*_\tau\|_{L^2(H^1)}^2
    + \|\bar{u}^*_\tau\|_{L^{2\tilde{q}_d/p_d}(L^{p_d})}^2 
    + \tau^{p_d/\tilde{q}_d}\|u_0\|_{L^{p_d}}^2\right)\\
    &\le c_4\left(\|\bar{u}^*_\tau\|_{L^2(H^1)}^2
    + c_3\|\bar{u}^*_\tau\|_{L^2(H^1)}^{2p_d/\tilde{q}_d}
    \|\bar{u}^*_\tau\|_{L^\infty(L^2)}^{2-2p_d/\tilde{q}_d}
    + \tau^{p_d/\tilde{q}_d}\|u_0\|_{L^{p_d}}^2\right),
\end{aligned}\]
\[\begin{aligned}
    &\quad\|\bar{h}_\tau\|_{L^{4/{p_d}}(L^2)}
    \le \left(\tau\sum^N_{k=1}
    \|\upre_k\|_{L^{\tilde{q}_d}}^{4/{p_d}}\|\upre_{k-1}\|_{L^{p_d}}^{4/{p_d}}\right)^{{p_d}/4}
    \le c_5\left(\tau\sum^N_{k=1}
    \|\upre_k\|_1^{4/{p_d}}\|\upre_{k-1}\|_{L^{p_d}}^{4/{p_d}}\right)^{{p_d}/4}\\
    &\le c_5\left(\|\bar{u}^*_\tau\|_{L^2(H^1)}^2
    + c_3\|\bar{u}^*_\tau\|_{L^2(H^1)}^{2p_d/\tilde{q}_d}
    \|\bar{u}^*_\tau\|_{L^\infty(L^2)}^{2-2p_d/\tilde{q}_d}
    + \tau^{p_d/\tilde{q}_d}\|u_0\|_{L^{p_d}}^2\right)
\end{aligned}\]
for constants $c_4$ and $c_5$. Hence, by (\ref{ineq_bdd}), 
the weak convergences (\ref{convgw}) and (\ref{convhw}) hold.
Moreover, since there exists a constant $c_6>0$ such that 
$|(g_k,v)|
\le \|g_k\|_{L^{\tilde{p}_d}}\|v\|_{L^{\tilde{q}_d}}
\le c_6\|g_k\|_{L^{\tilde{p}_d}}\|v\|_1$
for all $k=1,2,\ldots,N$ and $v\in\Hb^d$, we have
\[\begin{aligned}
    &\quad\left\|\frac{\partial \hat{u}_\tau}{\partial t}\right\|_{L^{4/{p_d}}(V^*)}
    = \biggl\{\int^T_0 \biggl(\sup_{0\ne v\in V}\frac{1}{\|v\|_1}
    \biggl|-a(\bar{u}^*_\tau(t),v)
    -(\bar{g}_\tau(t),v)-(\bar{h}_\tau(t),\nabla v)\\
    &\qquad\qquad\qquad\qquad+\langle f_\tau(t),v\rangle_H
    -\int_{\Gamma_2}p^b_\tau(t) v\cdot nds\biggr|\biggr)^{4/{p_d}}dt\biggr\}^{{p_d}/4}\\
    &= \left\{\int^T_0 \left(
    c_a\|\bar{u}^*_\tau(t)\|_1+c_6\|\bar{g}_\tau(t)\|_{L^{\tilde{p}_d}}
    +\|\bar{h}_\tau(t)\|_0+\|f_\tau(t)\|_{H^*}
    +\|p^b_\tau(t)\|_1\right)^{4/{p_d}}dt\right\}^{{p_d}/4}\\
    &\le\! T^{p_d/(2\tilde{q}_d)}(c_a\sqrt{c_1}
    \!+\! \|f\|_{L^2(H^*)}
    \!+\! \|p^b\|_{L^2(H^1)})
    \!+\! c_6\|\bar{g}_\tau\|_{L^{4/{p_d}}(L^{\tilde{p}_d})}
    \!+\! \|\bar{h}_\tau\|_{L^{4/{p_d}}(L^2)},\\
    &\quad\, \left\|\frac{\partial \hat{u}^\circ_\tau}{\partial t}
    -\frac{\partial \hat{u}_\tau}{\partial t}\right\|_{L^{4/{p_d}}(V^*)}    = \left\{\tau\sum^N_{k=1}\left(\sup_{0\ne v\in V}\frac{
        |(u^\circ_k-u^\circ_{k-1}-u_k+u_{k-1},v)|
    }{\tau\|v\|_1}
    \right)^{4/{p_d}}\right\}^{{p_d}/4}\\
    &\le 2\left(\tau\sum^N_{k=1}\sup_{0\ne v\in V}
    \frac{|(\nabla P_k,v)|^{4/{p_d}}}{\|v\|_1^{4/{p_d}}}
    + \tau\|\upre_1-u_0\|_{V^*}^{4/{p_d}} \right)^{{p_d}/4}\\
    &\le 2\left(\tau\sum^N_{k=1}\sup_{0\ne v\in V}
    \frac{1}{\|v\|_1^{4/{p_d}}}\left|\int_{\Gamma_2} p^b(t_k) v\cdot nds\right|^{4/{p_d}}
    \right)^{{p_d}/4}
    + 2\tau^{p_d/4}\|\upre_1-u_0\|_{V^*} \\
    &\le 2\left(\tau\sum^N_{k=1}\|p^b(t_k)\|_1^{4/{p_d}}\right)^{{p_d}/4}
    + 2c_2\tau^{p_d/4+1/2}
    \le 2T^{p_d/(2\tilde{q}_d)}(\|p^b\|_{L^2(H^1)} + c_2 T),
\end{aligned}\]
and $\|\frac{\partial \hat{u}_\tau}{\partial t}\|_{L^{4/{p_d}}(V^*)}$
and $\|\frac{\partial \hat{u}^\circ_\tau}{\partial t}\|_{L^{4/{p_d}}(V^*)}$
are also bounded. Hence, (\ref{convnhwvs}) holds.
Furthermore, $\|\hat{u}^\circ_\tau\|_{L^2(H^1)}$ is bounded:
by (\ref{ineq_bdd}) and (\ref{ineq_U0}),
\[\begin{aligned}
    \|\hat{u}^\circ_\tau\|_{L^2(H^1)}^2
    &= \sum^N_{k=1}\int^1_0
    \|(1-s)u^\circ_{k-1} + su^\circ_k\|_1^2 \tau ds\\
    &\le \sum^N_{k=1}\tau(\|u^\circ_{k-1}\|_1^2
    + \|u^\circ_k\|_1^2)\int^1_0\{(1-s)^2+s^2\} ds
    \le \frac{4}{3}\|\bar{u}^*_\tau\|_{L^2(H^1)}^2
    + \frac{2c_1}{3}
    \le 2c_1,
\end{aligned}\]
which implies the strong convergence (\ref{convshsl}) 
of $\hat{u}^\circ_\tau$ in $L^2(\Lo^d)$ from the Aubin--Lions lemma
\cite[Theorem II.5.16 (i)]{BF13}.
Since we have for all $t\in(t_{k-1},t_k), k=1,2,\ldots,N$,
\[\begin{aligned}
    \|\bar{u}^*_\tau(t)-\hat{u}^\circ_\tau(t)\|_0
    &= \left|\frac{t_k-t}{\tau}\right| \|u^\circ_k-u^\circ_{k-1}\|_0
    \le \|u^\circ_k-u_k\|_0+\tau\|D_\tau u_k\|_0+\|u_{k-1}-u^\circ_{k-1}\|_0,\\
    \|\bar{u}^*_\tau(t)-\hat{u}_\tau(t)\|_0
    &\le \|\bar{u}^*_\tau(t)-\bar{u}_\tau(t)\|_0
    +\|\bar{u}_\tau(t)-\hat{u}_\tau(t)\|_0
    \le \|\upre_k-u_k\|_0+\tau\|D_\tau u_k\|_0,
\end{aligned}\]
the functions $\bar{u}^*_{\tau_k}$, $\hat{u}^\circ_{\tau_k}$
and $\hat{u}_{\tau_k}$ possess a common limit function $u$,
and the strong convergences (\ref{convsbsl}) and (\ref{convnhsl}) hold:
by (\ref{ineq_bdd}) and (\ref{ineq_U0}),
\[\begin{aligned}
    \|\bar{u}^*_\tau-\hat{u}^\circ_\tau\|_{L^2(L^2)}
    &\le \left(\tau\sum_{k=1}^N(\|u^\circ_k-u_k\|_0
    + \tau\|D_\tau u_k\|_0 + \|u_{k-1}-u^\circ_{k-1}\|_0)^2\right)^{1/2}\\
    &\le 2\sqrt{3}\|\bar{u}^*_\tau-\bar{u}_\tau\|_{L^2(L^2)}
    \!+\! \sqrt{3}\tau\left\|\frac{\partial \hat{u}_\tau}{\partial t}\right\|_{L^2(L^2)}
    \!+\! \sqrt{3\tau} \|u_0-\upre_1\|_0
    \le 5\sqrt{3c_1 \tau},
\end{aligned}\]
\[\begin{aligned}
    \|\bar{u}^*_\tau-\hat{u}_\tau\|_{L^2(L^2)}
    &\le \sqrt{2}\|\bar{u}^*_\tau-\bar{u}_\tau\|_{L^2(L^2)}
    +\sqrt{2}\tau\left\|\frac{\partial \hat{u}_\tau}{\partial t}\right\|_{L^2(L^2)}
    \le 2\sqrt{2c_1\tau}.
\end{aligned}\]
It also holds that
\[
    \|\bar{u}^*_\tau-\hat{u}^\circ_\tau\|_{L^\infty(L^2)}
    \le \max_{k=1,2,\ldots,N}(\|u^\circ_k\|_0+\|u^\circ_{k-1}\|_0)
    \le 2\sqrt{c_1}.
\]
Since $\|\hat{u}^\circ_\tau\|_{L^\infty(L^2)}$ and 
$\|\frac{\partial \hat{u}^\circ_\tau}{\partial t}\|_{L^{4/p_d}(V^*)}$ 
are bounded, we obtain the strong convergence (\ref{convshsc})
of $\hat{u}^\circ_\tau$ in $C([0,T];V^*)$
\cite[Theorem II.5.16 (ii)]{BF13}.
In particular, $\hat{u}^\circ_\tau(0)$ converges to $u(0)$ in $V^*$.
On the other hand, by (\ref{ineq_Ustar}), $\hat{u}^\circ_\tau(0)=\upre_1$
converges to $u_0$ in $V^*$. Through the uniqueness of the limit in $V^*$,
we have indeed obtained that $u(0)=u_0$.

From (\ref{eq_ghver}) with $\eps:=\eps_k$, taking $k\rightarrow\infty$, 
it holds that for all $v\in V$ and $\theta\in C^\infty_0(0,T)$,
\[\begin{aligned}
    \int^T_0\left(\left\langle\frac{\partial u}{\partial t},
      \theta v\right\rangle_V
    + a(u,\theta v) + (g,\theta v) + (h,\nabla (\theta v))\right) dt
    =\int^T_0\left(\langle f,\theta v\rangle_H 
    - \int_{\Gamma_2}p^b\theta v\cdot nds\right)dt.
\end{aligned}\]
Next, we show that 
\begin{align}\label{eq_ghu}
    g = (\nabla u)^T u - u\dvg u,\qquad
    h = -u(u)^T.
\end{align}
We set $\bar{v}_\tau(t):=\upre_{k-1}$ for $t\in(t_{k-1},t_k],k=1,2,\ldots,N$.
Then it holds that 
\[\begin{aligned}
    \|\bar{v}_\tau-\bar{u}^*_\tau\|_{L^2(L^2)}
    &\!\le\! \left(\tau\sum_{k=1}^N(\|\upre_k-u_k\|_0
    + \tau\|D_\tau u_k\|_0 + \|u_{k-1}-\upre_{k-1}\|_0)^2\right)^{1/2}
    \!\le\! 3\sqrt{3c_1 \tau},
\end{aligned}\]
and hence it follows from (\ref{convsbsl}) that 
$\bar{v}_{\tau_k}\rightarrow u$ strongly in $L^2(\Lo^d)$ as $k\rightarrow\infty$. 
Since $\nabla\bar{u}^*_\tau\rightharpoonup \nabla u$ weakly in $L^2(\Lo^{d\times d})$
and $\dvg\bar{u}^*_\tau\rightharpoonup\dvg u$ weakly in $L^2(\Lo)$
as $k\rightarrow\infty$, we have
\[\begin{array}{rcll}
    \bar{g}_\tau 
    = (\nabla\bar{u}^*_\tau)^T \bar{v}_\tau 
    - \bar{v}_\tau\dvg\bar{u}^*_\tau
    &\rightharpoonup & (\nabla u)^T u - u\dvg u
    &\text{weakly in }L^1(L^1(\Omega)^d),\\[4pt]
    \bar{h}_\tau = -\bar{u}^*_\tau(\bar{v}_\tau)^T
    &\rightarrow & -u(u)^T
    &\text{strongly in }L^1(L^1(\Omega)^{d\times d})
\end{array}\]
as $k\rightarrow\infty$ (cf. \cite[Proposition II.2.12]{BF13}).
On the other hand, we also know (\ref{convgw}) and (\ref{convhw}).
The convergence in these spaces imply the convergence 
in the distributions sense, therefore (\ref{eq_ghu}) holds 
by the uniqueness of the limit in $\mathcal{D}'((0,T)\times\Omega)$.
Hence, it holds that for all $v\in V$ and $\theta\in C^\infty_0(0,T)$,
\[\begin{aligned}
    &\quad\int^T_0\left(\left\langle\frac{\partial u}{\partial t},v\right\rangle_V
    +a(u,v)+((\nabla u)u-u\dvg u,v) - (u(u)^T,\nabla v)\right)\theta dt\\
    &=\int^T_0\left(\langle f,v\rangle_H - \int_{\Gamma_2}p^bv\cdot nds\right)\theta dt,
\end{aligned}\]
which is equivalent to the following
\[
    \int^T_0\left(\left\langle\frac{\partial u}{\partial t},v\right\rangle_V
    + a(u,v) + d(u,u,v)\right)\theta dt
    =\int^T_0\left(\langle f,v\rangle_H 
    - \int_{\Gamma_2}p^bv\cdot nds\right)\theta dt.
\]
\qed

\end{document}